\documentclass[10pt]{amsart}
\usepackage{amscd}
\usepackage{amssymb}
\usepackage[all]{xy}
\date{7 August 2007}
\title{Rigid Dualizing Complexes over Commutative Rings}

\author{Amnon Yekutieli and James J. Zhang }
\address{A. Yekutieli: Department of  Mathematics 
Ben Gurion University, 
Be'er Sheva 84105, 
Israel}
\email{amyekut@math.bgu.ac.il}
\address{J.J. Zhang: Department of Mathematics, Box 354350,
University of Washington, Seattle, Washington 98195, USA}
\email{zhang@math.washington.edu}
\thanks{{\em Mathematics Subject Classification} 2000.
Primary: 14F05; Secondary: 14B25, 14F10, 13D07, 18G10, 16E45.}
\keywords{commutative rings, derived categories,
dualizing complexes, rigid complexes.}
\thanks{This research was supported by the US-Israel Binational
Science Foundation. The second author was partially supported by 
the US National Science Foundation.}

\newtheorem{thm}[equation]{Theorem}
\newtheorem{cor}[equation]{Corollary}
\newtheorem{prop}[equation]{Proposition}
\newtheorem{lem}[equation]{Lemma}
\theoremstyle{definition}
\newtheorem{dfn}[equation]{Definition}
\newtheorem{rem}[equation]{Remark}
\newtheorem{exa}[equation]{Example}

\numberwithin{equation}{section}

\newcommand{\iso}{\xrightarrow{\simeq}}

\newcommand{\xar}{\xrightarrow}

\newcommand{\opn}{\operatorname}
\newcommand{\cat}[1]{\operatorname{\mathsf{#1}}}

\newcommand{\rmitem}[1]{\item[\text{\textup{(#1)}}]}
\newcommand{\mfrak}[1]{\mathfrak{#1}}
\newcommand{\mcal}[1]{\mathcal{#1}}
\newcommand{\msf}[1]{\mathsf{#1}}

\newcommand{\mrm}[1]{\mathrm{#1}}
\newcommand{\mbb}[1]{\mathbb{#1}}

\newcommand{\tup}[1]{\textup{#1}}
\newcommand{\bsym}[1]{\boldsymbol{#1}}
\newcommand{\boplus}{\bigoplus\nolimits}

\newcommand{\til}[1]{\tilde{#1}}
\newcommand{\what}[1]{\widehat{#1}}

\renewcommand{\k}{\Bbbk}
\newcommand{\K}{\mbb{K} \hspace{0.05em}}
\renewcommand{\d}{\mathrm{d}}

\newcommand{\cmnt}[1] 
{\mbox{} \newline 
\marginpar{\sf $\Lleftarrow$ comment } 
\textsf{[#1]} \newline}

\begin{document}

\begin{abstract}
In this paper we present a new approach to Grothendieck duality 
over commutative rings. 
Our approach is based on the idea of rigid dualizing 
complexes, which was introduced by Van den Bergh in the context 
of noncommutative algebraic geometry.
The method of rigidity was modified to work over 
general commutative base rings in our paper \cite{YZ5}.
In the present paper we obtain many of the 
important local features of Grothendieck duality, 
yet manage to avoid lengthy and difficult compatibility 
verifications. Our results apply to essentially finite type 
algebras over a regular noetherian finite dimensional base ring, 
and hence are suitable for arithmetic rings. In the sequel paper 
\cite{Ye4} these results will be used to construct and study rigid 
dualizing complexes on schemes.
\end{abstract}

\maketitle

\setcounter{section}{-1}
\section{Introduction}

Grothendieck duality for schemes was introduced in the book 
``Residues and Duality'' \cite{RD} by R. Hartshorne. This duality 
theory has applications in various areas of algebraic geometry, 
including moduli spaces, resolution of singularities, arithmetic 
geometry, enumerative geometry and more. 

In the forty years since the publication of \cite{RD} a number of 
related papers appeared in the literature. Some of these papers 
provided elaborations on, or more explicit versions of Grothendieck 
duality (e.g.\ \cite{Kl}, \cite{Li1}, \cite{HK}, \cite{Ye2}, 
\cite{Ye3}, \cite{Sa}). Other papers contained alternative 
approaches (e.g.\ \cite[Appendix]{RD}, \cite{Ve} and \cite{Ne}). 
The recent book \cite{Co} is a complement to \cite{RD} that 
fills gaps in the proofs, and also contains the first proof of 
the Base Change Theorem. A noncommutative version of Grothendieck 
duality was developed in \cite{Ye1}, which has applications in 
algebra (e.g.\ \cite{EG}) and even in mathematical physics (e.g.\ 
\cite{KKO}). Other papers sought to extend the scope of 
Grothendieck duality to formal schemes (e.g.\ \cite{AJL} and 
\cite{LNS}) or to differential graded algebras (see \cite{FIJ}).

One of the fascinating features of Grothendieck duality
is the complicated interplay between its local and global 
components. Another feature of this theory (in some sense 
parallel to the first) is the gap between formal categorical 
statements and their concrete realizations. 
Much of the effort in studying Grothendieck duality was aimed at 
clarifying the local-global interplay, and at bridging the 
above-mentioned gap. 

In this paper we present a new approach to Grothendieck duality 
for commutative rings (i.e.\ affine schemes). The sequel
\cite{Ye4} will treat schemes in general (including duality for 
proper morphisms). The key idea in our approach is the use of 
{\em rigid dualizing 
complexes}. This notion was introduced by Van den Bergh 
\cite{VdB} in the context of noncommutative algebraic geometry, 
and was developed further, in the noncommutative direction,
in our papers \cite{YZ1, YZ2, YZ3, YZ4}. 
In the paper \cite{YZ5} we worked out the fundamental properties 
of rigid complexes over commutative rings relative to an 
arbitrary commutative base ring (as opposed to a base field). 
Attaching a rigid structure to a dualizing complex eliminates all 
nontrivial automorphisms. Moreover, rigid dualizing complexes 
admit several useful operations, such as localization and traces.

The general concept of ``rigidity'' is familiar in other areas 
of algebraic geometry (e.g.\ level structures on elliptic curves, 
or marked points on higher genus curves). Actually, 
in Grothendieck's original treatment \cite{RD} duality 
(local and global) itself was used as a sort of ``rigid structure''
on dualizing complexes, but 
this was very cumbersome (amounting to big commutative 
diagrams) and hard to employ. On the other hand, Van den Bergh's 
rigidity is very neat, and enjoys remarkable functorial properties. 

The background material we need in this paper 
is standard commutative algebra, 
the theory of derived categories, and the theory of rigid 
complexes over commutative rings from \cite{YZ5}. All of that is 
reviewed is Section 1 of the paper, for the convenience of the 
reader. We also need a few isolated 
results on dualizing complexes from \cite{RD}. 

Let us explain what are rigid dualizing complexes and how they
are used in our paper. Fix for the rest of the introduction a 
finite dimensional, regular, noetherian, commutative base ring 
$\K$ (e.g.\ a field, or the ring of integers). 
Recall that an essentially finite type commutative
$\K$-algebra $A$ is by definition a localization 
(with respect to a multiplicatively closed subset)
of some finitely 
generated $\K$-algebra. Note that the rings 
$\Gamma(U, \mcal{O}_{X})$ and $\mcal{O}_{X, x}$, where $X$ is a 
finite type $\K$-scheme, $U \subset X$ is an affine open set and 
$x \in X$ is a point, are both essentially finite type $\K$-algebras.
We denote by $\cat{EFTAlg} / \K$ the category of essentially
finite type commutative
$\K$-algebras, and by default we will stay within this 
category. 

We shall use notation such as $f^* : A \to B$ for a 
homomorphism of algebras, corresponding to a morphism of schemes
$f : \opn{Spec} B \to \opn{Spec} A$.
This convention, although perhaps awkward at first sight, fits 
better with the usual notation for associated functors.
Thus there are functors
$f_* : \cat{Mod} B \to \cat{Mod} A$
(restriction of scalars, which is push-forward geometrically) and
$f^* : \cat{Mod} A \to \cat{Mod} B$
(extension of scalars, i.e.\ $B \otimes_A -$,
which is pull-back geometrically).
For composable homomorphisms 
$A \xar{f^*} B \xar{g^*} C$ we sometimes write
$(f \circ g)^*$ instead of $g^* \circ f^*$. 

For a $\K$-algebra $A$ the derived category of  
complexes of $A$-modules is denoted by
$\msf{D}(\cat{Mod} A)$, with the usual modifiers (e.g.\ 
$\msf{D}^{\mrm{b}}_{\mrm{f}}(\cat{Mod} A)$ is the full subcategory 
of bounded complexes with finitely generated cohomologies).

Given a complex 
$M \in \msf{D}(\cat{Mod} A)$
we define its {\em square}
$\opn{Sq}_{A / \K} M \in \msf{D}(\cat{Mod} A)$.
The functor $\opn{Sq}_{A / \K}$ from the category 
$\msf{D}(\cat{Mod} A)$ to itself is quadratic, in the sense that 
given a morphism 
$\phi : M \to N$ in $\msf{D}(\cat{Mod} A)$
and an element $a \in  A$, one has
$\opn{Sq}_{A / \K}(a \phi) = a^2 \opn{Sq}_{A / \K}(\phi)$.

A {\em rigidifying isomorphism} for $M$ is an isomorphism
$\rho : M \iso \opn{Sq}_{A / \K} M$
in $\msf{D}(\cat{Mod} A)$. If 
$M \in \msf{D}^{\mrm{b}}_{\mrm{f}}(\cat{Mod} A)$
then the pair $(M, \rho)$ 
is called a {\em rigid complex over $A$ relative to $\K$}. 
Suppose $(M, \rho_M)$ and $(N, \rho_N)$ are two rigid complexes. 
A {\em rigid morphism} 
$\phi : (M, \rho_M) \to (N, \rho_N)$ is 
a morphism $\phi : M \to N$ in $\msf{D}^{}(\cat{Mod} A)$
such that 
$\rho_N \circ \phi = \opn{Sq}_{A / \K}(\phi) \circ \rho_M$.
Observe that if $(M, \rho_M)$ is a rigid complex such that
$\opn{RHom}_{A}(M, M) = A$, and 
$\phi : (M, \rho_M) \to (M, \rho_M)$ is a rigid isomorphism, then 
$\phi$ is multiplication by some invertible element $a \in A$ 
satisfying $a = a^2$; and therefore $a = 1$. We conclude that 
{\em the identity is the only rigid automorphism of 
$(M, \rho_M)$.}

Let $B$ be another $\K$-algebra, 
and let $f^* : A \to B$ be a finite $\K$-algebra
homomorphism. Define
\[ f^{\flat} M := \opn{RHom}_{A}(B, M) \in 
\msf{D}^{+}_{\mrm{f}}(\cat{Mod} B) . \]
If $f^{\flat} M$ has bounded cohomology
then there is an induced rigidifying isomorphism 
$f^{\flat}(\rho_M) : f^{\flat} M \iso 
\opn{Sq}_{B / \K} f^{\flat} M$.
We write
$f^{\flat}(M, \rho_M) :=  
(f^{\flat} M, f^{\flat}(\rho_M))$.

Next we consider {\em essentially smooth} 
homomorphisms. By definition 
$f^* : A \to B$ is essentially smooth if it is essentially 
finite type and formally smooth. Then $B$ is flat over $A$, and
the module of differentials 
$\Omega^1_{B / A}$ is a finitely generated projective 
$B$-module. The rank of $\Omega^1_{B / A}$ might vary on 
$\opn{Spec} B$; but if it has constant rank $n$ then we say 
$f^*$ is essentially smooth of relative dimension $n$. 
An essentially smooth homomorphism 
of relative dimension $0$ is called
an {\em essentially \'etale} homomorphism.

Let $f^* : A \to B$ be an essentially smooth homomorphism,
let $B = \prod B_i$ be the decomposition of $\opn{Spec} B$ into 
connected components, and for any $i$ let
$n_i := \opn{rank}_{B_i} \Omega^1_{B_i / A}$. 
We then define
\[ f^{\sharp} M := \boplus_i \Omega^{n_i}_{B_i / A}[n_i] 
\otimes_A M \in \msf{D}^{\mrm{b}}_{\mrm{f}}(\cat{Mod} B) . \]
There is an induced rigidifying isomorphism
$f^{\sharp}(\rho_M) : f^{\sharp} M \iso 
\opn{Sq}_{B / \K} f^{\sharp} M$, 
and thus a new rigid complex
$f^{\sharp}(M, \rho_M) := (f^{\sharp} M, f^{\sharp}(\rho_M))$.

Now let's consider dualizing complexes. 
Recall that a complex 
$R \in \msf{D}^{\mrm{b}}_{\mrm{f}}(\cat{Mod} A)$ is dualizing if 
it has finite injective dimension over $A$, 
and if the canonical morphism 
$A \to \opn{RHom}_{A}(R, R)$ is an isomorphism.
A {\em rigid dualizing complex over $A$ relative to $\K$} is a 
rigid complex $(R, \rho)$ such that $R$ is dualizing. 

Here is the main result of our paper.

\begin{thm} \label{thm0.1}
Let $\K$ be a regular finite dimensional noetherian ring, and
let $A$ be an essentially finite type $\K$-algebra.
\begin{enumerate}
\item The algebra $A$ has a rigid dualizing complex 
$(R_A, \rho_A)$, which is unique up to a unique rigid isomorphism.
\item Given a finite homomorphism $f^* : A \to B$, there is a 
unique rigid isomorphism
$f^{\flat}(R_A, \rho_A) \iso (R_B, \rho_B)$.
\item Given an essentially smooth homomorphism $f^* : A \to B$ , 
there is a unique rigid isomorphism
$f^{\sharp}(R_A, \rho_A) \iso (R_B, \rho_B)$.
\end{enumerate}
\end{thm}

This theorem is repeated as Theorem \ref{thm2.5}
in the body of the paper. 

The next result is one we find quite surprising. Its significance 
is not yet understood. 

\begin{thm}
Let $\K$ be a regular finite dimensional noetherian ring, and
let $A$ be an essentially finite type $\K$-algebra.
Assume $\opn{Spec} A$ is connected and nonempty. 
Then, up to isomorphism, 
the only nonzero rigid complex over $A$ relative to $\K$ is
the rigid dualizing complex $(R_A, \rho_A)$.
\end{thm}

This theorem is repeated as Theorem \ref{thm2.4}. 

Let $A$ be some $\K$-algebra. Using the rigid dualizing complex 
$R_A$ we define the {\em auto-duality functor} 
$\mrm{D}_A := \opn{RHom}_A(-, R_A)$
of $\msf{D}_{\mrm{f}}(\cat{Mod} A)$.
Due to Theorem \ref{thm0.1}(1) this functor
is independent of the rigid dualizing complex $R_A$ chosen.
Given a $\K$-algebra homomorphism $f^* : A \to B$ 
we define the {\em twisted inverse image functor}
\[ f^! : \msf{D}^{+}_{\mrm{f}}(\cat{Mod} A) \to
\msf{D}^{+}_{\mrm{f}}(\cat{Mod} B)  \]
as follows. If $A=B$ and $f^* = \bsym{1}_A$ (the identity 
automorphism) then
$f^! := \bsym{1}_{\msf{D}^{+}_{\mrm{f}}(\cat{Mod} A)}$
(the identity functor). Otherwise we define
$f^! := \mrm{D}_B \, \mrm{L} f^* \, \mrm{D}_A$.

As explained in Corollary \ref{cor3.1},
the base ring $\K$ can sometimes be ``factored out'' of 
the construction of the twisted inverse image functor $f^!$. 

Suppose we are given two homomorphisms 
$A \xar{f^*} B \xar{g^*} C$ in $\cat{EFTAlg} / \K$.
Then there is an obvious isomorphism 
\[ \phi_{f,g} : (f \circ g)^! \iso g^! f^! \]
of functors 
$\msf{D}^{+}_{\mrm{f}}(\cat{Mod} A) \to
\msf{D}^{+}_{\mrm{f}}(\cat{Mod} C)$,
coming from the adjunction isomorphism
$\bsym{1} \iso \mrm{D}_B\, \mrm{D}_B$ on 
$\msf{D}^{+}_{\mrm{f}}(\cat{Mod} B)$. 
For three homomorphisms 
$A \xar{f^*} B \xar{g^*} C \xar{h^*} D$ in 
$\cat{EFTAlg} / \K$ the isomorphisms $\phi_{-,-}$ satisfy
the compatibility condition 
\[ \phi_{g, h} \circ \phi_{f, g \circ h}  
= \phi_{f, g} \circ \phi_{f \circ g, h} 
: (f \circ g \circ h)^! \iso h^! g^! f^!  \]
(see Proposition \ref{prop3.1}). 
This means that the assignment 
$f^* \mapsto f^!$ is the $1$-component 
of a $2$-functor $\cat{EFTAlg} / \K \to \cat{Cat}$, 
whose $0$-component is 
$A \mapsto \msf{D}^{+}_{\mrm{f}}(\cat{Mod} A)$. 
Here $\cat{Cat}$ denotes the $2$-category of all categories.
(The notion of $2$-functor is recalled in Section 
\ref{sec.2-func}.)
 
The next theorem (which is a consequence of 
Theorem \ref{thm0.1}) describes the variance properties of the
twisted inverse image $2$-functor $f^* \mapsto f^{!}$.

\begin{thm} \label{thm0.9}
Let $f^* : A \to B$ be a homomorphism in 
$\cat{EFTAlg} / \K$. 
\begin{enumerate}
\item If $f^*$ is finite, then there is an isomorphism 
\[ \psi_f^{\flat} : f^{\flat} \iso f^! \]
of functors 
$\msf{D}^{+}_{\mrm{f}}(\cat{Mod} A) \to 
\msf{D}^{+}_{\mrm{f}}(\cat{Mod} B)$. 
These isomorphisms are $2$-fun\-ctorial for finite homomorphisms. 
\item If $f^*$ is essentially smooth, then there is an isomorphism 
\[ \psi_f^{\sharp} : f^{\sharp} \iso f^!  \]
of functors $\msf{D}^{+}_{\mrm{f}}(\cat{Mod} A) \to 
\msf{D}^{+}_{\mrm{f}}(\cat{Mod} B)$. 
These isomorphisms are $2$-fun\-ctorial for essentially 
smooth homomorphisms. 
\end{enumerate}
\end{thm}

For more detailed statements and proofs 
see Theorems \ref{thm3.2} and \ref{thm3.3}.

In the situation of a finite homomorphism $f^* : A \to B$, part 
(1) of the theorem gives rise to the {\em functorial trace map}
\[ \opn{Tr}_f : f_* f^! \to \bsym{1} . \]
This is a nondegenerate morphism of functors from 
$\msf{D}^{+}_{\mrm{f}}(\cat{Mod} A)$ to itself. See Proposition 
\ref{prop.funct-traces.1} for details. 

If $f^* : A \to B$ is essentially \'etale, then from part 
(2) of the theorem we get the {\em functorial localization map}
\[ \opn{q}_f : \bsym{1} \to f_* f^!  , \]
which is a nondegenerate morphism of functors from 
$\msf{D}^{+}_{\mrm{f}}(\cat{Mod} A)$ to itself. See Proposition 
\ref{prop4.1} for details. 
The relation between functorial localization maps and 
functorial trace maps is explained in Proposition 
\ref{prop.funct-traces.3}.

Here is an application to differential forms, which is a 
corollary to Theorem \ref{thm0.9}.

\begin{cor} \label{cor0.1}
Suppose $A \to B \to C$ are homomorphisms in 
$\cat{EFTAlg} / \K$, with $A \to B$ and $A \to C$ essentially
smooth of relative dimension $n$, and $B \to C$ finite. 
Then there is a nondegenerate trace map
\[ \opn{Tr}_{C / B / A} : \Omega^n_{C/A} \to \Omega^n_{B/A} . \]
The trace maps $\opn{Tr}_{- / - / A}$ are
functorial for such finite homomorphisms
$B \to C$, and commute with the localization maps
for a localization homomorphism $B \to B'$. 
\end{cor}

This is restated (in more detail) as Theorem \ref{thm5.1} and 
Proposition \ref{prop5.3}. In some cases we can compute 
$\opn{Tr}_{C / B / A}$ ; see 
Propositions \ref{prop5.1} and \ref{prop5.2}. 

For a finite 
flat homomorphism $A \to B$ let $\opn{tr}_{B / A} : B \to A$ be the 
usual trace map, i.e.\ $\opn{tr}_{B / A}(b) \in A$ 
is the trace of the
operator $b$ acting of the locally free $A$-module $B$. 

Here is another result relating traces and localization. 
It is not a consequence of Theorem \ref{thm0.9}, but instead 
relies on the fact that for an \'etale homomorphism $A \to B$
one has a canonical ring isomorphism
$B \otimes_A B \cong B \times B'$, where $B'$ is the kernel of 
the multiplication map $B \otimes_A B \to B$. 
This is interpreted in terms of rigidity.

\begin{thm} \label{thm0.6}
Suppose $f^* : A \to B$ is a finite \'etale homomorphism in
$\cat{EFTAlg} / \K$, so the localization map
$\opn{q}_f : A \to f^! A$ induces an isomorphism 
$\bsym{1} \otimes \opn{q}_f : B \iso f^! A$. 
Then the diagram
\[ \UseTips \xymatrix @C=5ex @R=3.5ex {
B 
\ar[rr]^{\bsym{1} \otimes \opn{q}_f} 
\ar[dr]_{\opn{tr}_{B/A}} 
& & f^! A
\ar[dl]^{\opn{Tr}_{f; A}} \\
& A
} \]
is commutative.
\end{thm}

This result is repeated (in slightly more general form) 
as Theorem \ref{thm3.1}.

To conclude the introduction let us mention how our twisted 
inverse image $2$-functor compares with the original constructions 
in \cite{RD}. We shall restrict attention to the category 
$\cat{FTAlg} / \K$ of finite type $\K$-algebras.
(This is of course the opposite of the category of finite type 
affine $\K$-schemes.) Given a homomorphism
$f^* : A \to B$ in $\cat{FTAlg} / \K$ let us denote by
\[ f^{! (\mrm{G})} : \msf{D}^{+}_{\mrm{f}}(\cat{Mod} A) \to
\msf{D}^{+}_{\mrm{f}}(\cat{Mod} B) \]
the twisted inverse image from \cite{RD}.
In particular, for an algebra $A$, with structural homomorphism
$\pi_A^* : \K \to A$, we obtain the complex
$R^{(\mrm{G})}_A := \pi_A^{! (\mrm{G})} \K \in
\msf{D}^{+}_{\mrm{f}}(\cat{Mod} A)$,
which is known to be dualizing. In Theorem \ref{thm3.4}
we show that for any $A \in \cat{FTAlg} / \K$ there is an isomorphism
$R^{(\mrm{G})}_A \cong R_A$ in $\msf{D}(\cat{Mod} A)$.
This implies that there is an isomorphism 
$f^! \cong f^{! (\mrm{G})}$
of $2$-functors
$\cat{FTAlg} / \K \to \cat{Cat}$.
In general we do not know an easy way to make the isomorphisms
$R^{(\mrm{G})}_A \cong R_A$ canonical; but see Remark 
\ref{rem3.11}.

\medskip \noindent
\textbf{Acknowledgments.}
The authors wish to Brian Conrad, Reinhold H\"ubl, Steven Kleiman,
Joseph Lipman, Amnon Neeman, Paramathanath Sastry 
and Michel Van den Bergh 
for useful discussions and valuable suggestions.

\tableofcontents

\section{Review of Rigid Complexes}
\label{sec.rev}

In this section we recall definitions and results from the paper 
\cite{YZ5}. 

Throughout the paper all rings and algebras are assumed to be 
commutative by default.
Given two rings $A$ and $B$ we denote a 
homomorphism between them by an expression such as
$f^* : A \to B$. This of course signifies that the corresponding 
morphism of schemes is $f : \opn{Spec} B \to \opn{Spec} A$. 
The benefit of this notation is that it is compatible 
with the customary notation for various related functors, such as
$f_* : \cat{Mod} B \to \cat{Mod} A$ (restriction of scalars, 
which is direct image for schemes), and
$f^* :\cat{Mod} A \to \cat{Mod} B$ (base change, i.e.\ 
$f^* M = B \otimes_A M$, which is inverse image for schemes). 

Given a ring $A$ we denote by $\msf{D}(\cat{Mod} A)$ the derived 
category of complexes of $A$-modules. If $A$ is noetherian then we 
denote by $\msf{D}^{\mrm{b}}_{\mrm{f}}(\cat{Mod} A)$ the full 
subcategory consisting of bounded complexes with finitely 
generated cohomology modules. 

In \cite[Section 2]{YZ5} we introduced the {\em squaring 
operation}. Let $A$ be a ring, $B$ an $A$-algebra and 
$M \in \msf{D}(\cat{Mod} B)$. The square of $M$ over $B$ relative 
to $A$ is a complex 
$\opn{Sq}_{B / A} M \in \msf{D}(\cat{Mod} B)$.
In case $B$ is flat over $A$ one has the simple formula
\[ \opn{Sq}_{B / A} M = \opn{RHom}_{B \otimes_A B}
(B,M \otimes^{\mrm{L}}_{A} M) . \]
But in general it is necessary to replace 
$B \otimes_A B$ with $\til{B} \otimes_A \til{B}$
in the formula defining $\opn{Sq}_{B / A} M$,
where $\til{B}$ is a suitable differential graded $A$-algebra,
quasi-isomorphic to $B$. 

Suppose $C$ is another $A$-algebra, $f^* : B \to C$ is an 
$A$-algebra homomorphism,  
$N \in \msf{D}(\cat{Mod} C)$, and $\phi : N \to M$ is a morphism in 
$\msf{D}(\cat{Mod} B)$. (Strictly speaking, $\phi$ is a morphism
$f_* N \to M$). Then there is an induced morphism
\[ \opn{Sq}_{f^* / A}(\phi) : \opn{Sq}_{C / A} N \to
\opn{Sq}_{B / A} M \]
in $\msf{D}(\cat{Mod} B)$. The formation of 
$\opn{Sq}_{f^* / A}(\phi)$ is functorial in $f^*$ and in $\phi$. 

Specializing to the case $C = B$ and $f^* = \bsym{1}_B$ (the 
identity homomorphism), we obtain a functor
\[ \opn{Sq}_{B / A} := \opn{Sq}_{\bsym{1}_B / A} : 
\msf{D}(\cat{Mod} B) \to \msf{D}(\cat{Mod} B) . \]
This functor is quadratic, in the sense that for any
$\phi \in \opn{Hom}_{\msf{D}(\cat{Mod} B)}(M, N)$
and $b \in B$ one has
\[ \opn{Sq}_{B / A}(b \phi) = b^2 \opn{Sq}_{B / A}(\phi) . \]

The next definition is a variant of the original definition of 
Van den Bergh \cite{VdB}.

\begin{dfn} \label{dfn1.1}
Let $A$ be a ring and $B$ a noetherian $A$-algebra. A {\em rigid 
complex over $B$ relative to $A$} is a pair $(M, \rho)$, where:
\begin{enumerate}
\item $M$ is a complex in 
$\msf{D}^{\mrm{b}}_{\mrm{f}}(\cat{Mod} B)$ which has finite flat 
dimension over $A$.
\item $\rho$ is an isomorphism
\[ \rho : M \iso \opn{Sq}_{B / A} M \]
in $\msf{D}(\cat{Mod} B)$, called a {\em rigidifying isomorphism}.
\end{enumerate}
\end{dfn}

\begin{dfn} \label{dfn.rev.1}
Let $A$ be a ring, let $B$ and $C$ be noetherian $A$-algebras, let
$f^* : B \to C$ be an $A$-algebra homomorphism, let
$(M, \rho)$ be a rigid complex over $B$ relative to $A$, and let
$(N, \sigma)$ be a rigid complex over $C$ relative to $A$.
A {\em rigid trace morphism relative to $A$}
\[ \phi :   (N, \sigma) \to (M, \rho) \]
is a morphism $\phi : N \to M$ in $\msf{D}(\cat{Mod} B)$, such that 
the diagram
\[ \UseTips \xymatrix @C=5ex @R=5ex {
N 
\ar[r]^(0.35){\sigma}
\ar[d]_{\phi}
& \opn{Sq}_{C / A} N
\ar[d]^{\opn{Sq}_{f^* / A}(\phi)} \\
M
\ar[r]^(0.35){\rho}
& \opn{Sq}_{B / A} M
} \]
is commutative.
\end{dfn}

Clearly the composition of two rigid trace morphisms is again a 
rigid trace morphism.

Specializing Definition \ref{dfn.rev.1} 
to the case $C = B$ and $f^* = \bsym{1}_B$, 
we call such a morphism $\phi$ 
a {\em rigid morphism over $B$ relative to $A$}. Let 
$\msf{D}^{\mrm{b}}_{\mrm{f}}(\cat{Mod} B)_{\mrm{rig} / A}$
be the category whose objects are the rigid complexes 
over $B$ relative to $A$, and whose morphisms are the rigid 
morphisms. 

Recall that a ring homomorphism $f^* : A \to B$ is called 
{\em finite} if $B$ is a finitely generated $A$-module. 
Given a finite homomorphism $f^* : A \to B$ we define a functor
\[ f^{\flat} : \msf{D}(\cat{Mod} A) \to \msf{D}(\cat{Mod} B) \]
by
\[ f^{\flat} M := \opn{RHom}_A(B, M) . \]
For any $M \in \msf{D}(\cat{Mod} A)$ there is a morphism
\begin{equation} \label{eqn1.1}
\opn{Tr}^{\flat}_{f; M} : f_* f^{\flat} M \to M
\end{equation}
called the {\em trace map}, which is induced from the homomorphism
$\phi \mapsto \phi(1)$ for $\phi \in \opn{Hom}_A(B, M)$.
In this way we get a morphism
$\opn{Tr}^{\flat}_{f} : f_* f^{\flat} \to \bsym{1}$
of functors from $\msf{D}(\cat{Mod} A)$ to itself.
Note that for any $N \in \msf{D}(\cat{Mod} B)$ the homomorphism 
of $B$-modules
\[ \opn{Hom}_{\msf{D}(\cat{Mod} B)}(N, f^{\flat} M) \to
\opn{Hom}_{\msf{D}(\cat{Mod} A)}(N, M), \,
\psi \mapsto \opn{Tr}^{\flat}_{f; M} \circ \, \psi \]
is bijective (this is the adjunction isomorphism).

\begin{dfn} \label{dfn1.2}
Let $f^* : A \to B$ be a ring homomorphism, 
$M \in \msf{D}(\cat{Mod} A)$ 
and $N \in \msf{D}(\cat{Mod} B)$. A morphism
$\phi : N \to M$ in $\msf{D}(\cat{Mod} A)$ is called a 
{\em nondegenerate trace morphism} if the corresponding
morphism $N \to f^{\flat} M$ in 
$\msf{D}(\cat{Mod} B)$ is an isomorphism.
\end{dfn}

Here's the first result about rigid complexes, which explains 
their name.

\begin{thm}\label{thm.rev.1}
Let $A$ be a ring, let $f^* : B \to C$ be a 
homomorphism between noetherian $A$-algebras, let 
$(M, \rho) \in 
\msf{D}^{\mrm{b}}_{\mrm{f}}(\cat{Mod} B)_{\mrm{rig} / A}$
and let
$(N, \sigma) \in 
\msf{D}^{\mrm{b}}_{\mrm{f}}(\cat{Mod} C)_{\mrm{rig} / A}$.
Assume the canonical ring homomorphism 
$C \to \opn{End}_{\msf{D}(\cat{Mod} C)}(N)$
is bijective. Then there is at most one nondegenerate rigid trace 
morphism
$(N, \sigma) \to (M, \rho)$
over $B$ relative to $A$.
\end{thm}

\begin{proof}
This is a slight generalization of \cite[Theorem 0.2]{YZ5}.
Suppose we are given two nondegenerate rigid trace 
morphisms
$\phi, \phi' : (N, \sigma) \to (M, \rho)$. 
Since $\phi$ is nondegenerate it follows that 
$\opn{Hom}_{\msf{D}(\cat{Mod} B)}(N, M)$ is a free $C$-module with 
basis $\phi$. So $\phi' = c \phi$ for some unique element
$c \in C$. Because $\phi'$ is nondegenerate too we see that 
$c$ is an invertible element.
Next, since both $\phi$ and $\phi'$ are rigid, 
\cite[Corollary 2.7]{YZ5} says that $c^2 = c$. Thus $c = 1$
and $\phi' = \phi$.
\end{proof}

\begin{cor}[{\cite[Theorem 0.2]{YZ5}}] 
Taking $B = C$, $f^* = \bsym{1}_B$ and $(N, \sigma) = (M, \rho)$ 
in Theorem \tup{\ref{thm.rev.1}}, 
we see that the only automorphism of
$(M, \rho)$ in 
$\msf{D}^{\mrm{b}}_{\mrm{f}}(\cat{Mod} B)_{\mrm{rig} / A}$
is the identity $\bsym{1}_M$. 
\end{cor}

Here are a few results about pullbacks of rigid complexes.

\begin{thm}[{\cite[Theorem 5.3]{YZ5}}] \label{thm1.2}
Let $A$ be a noetherian ring, let $B, C$ be essentially finite 
type $A$-algebras, let $f^* : B \to C$ be a finite $A$-algebra
homomorphism, and let $(M, \rho) \in 
\msf{D}^{\mrm{b}}_{\mrm{f}}(\cat{Mod} B)_{\mrm{rig} / A}$.
Assume $f^{\flat} M$ has finite flat 
dimension over $A$.
\begin{enumerate}
\item The complex $f^{\flat} M$ has an induced rigidifying 
isomorphism 
\[ f^{\flat}(\rho) : f^{\flat} M  \iso \opn{Sq}_{C / A}
f^{\flat} M . \]
\item The rigid complex 
\[ f^{\flat}(M, \rho) := 
\bigl( f^{\flat} M, f^{\flat}(\rho) \bigr) \in
\msf{D}^{\mrm{b}}_{\mrm{f}}(\cat{Mod} C)_{\mrm{rig} / A} \]
depends functorially on $(M, \rho)$ and on $f^*$. 
\item Assume moreover that 
$\opn{Hom}_{\msf{D}(\cat{Mod} C)}(f^{\flat} M, f^{\flat} M) = C$. 
Then $\opn{Tr}^{\flat}_{f; M}$ is the unique nondegenerate rigid 
trace morphism
$f^{\flat} (M, \rho_M) \to (M, \rho_M)$
over $B$ relative to $A$.
\end{enumerate}
\end{thm}

Let $A$ be a noetherian ring. 
Recall that an $A$-algebra $B$ is called 
formally smooth (resp.\ formally
\'etale) if it has the lifting property (resp.\ the unique lifting 
property) for infinitesimal extensions. The $A$-algebra $B$ 
is called smooth (resp.\ \'etale) if it is finitely generated and
formally smooth (resp.\ formally \'etale).

In \cite[Section 3]{YZ5} 
we introduced a slightly more general kind of 
ring homomorphism than a smooth homomorphism. Again $A$ is 
noetherian. Recall that an $A$-algebra $B$ is called 
essentially finite type if it is a localization of some finitely 
generated $A$-algebra. We say that $B$ is {\em essentially smooth}
(resp.\ {\em essentially \'etale}) over $A$ if it is essentially 
finite type and formally smooth (resp.\ formally \'etale).
The composition of two 
essentially smooth homomorphisms is essentially smooth.
If $A \to B$ is essentially smooth then $B$ is flat over $A$, and
$\Omega^1_{B / A}$ is a finitely generated projective $B$-module.

Let $A$ be a noetherian ring and $f ^* : A \to B$ an 
essentially smooth homomorphism. Let
$\opn{Spec} B = \coprod_i \opn{Spec} B_i$
be the decomposition into connected components, and for every $i$ 
let $n_i$ be the rank of $\Omega^1_{B_i / A}$. We define a functor
\[ f^{\sharp} : \msf{D}(\cat{Mod} A) \to \msf{D}(\cat{Mod} B) \]
by
\[ f^{\sharp} M := \bigoplus_i \, 
\Omega^{n_i}_{B_i / A}[n_i] \otimes_A M . \]

\begin{thm}[{\cite[Theorem 6.3]{YZ5}}] \label{thm1.3}
Let $A$ be a noetherian ring, let $B, C$ be essentially finite 
type $A$-algebras, let $f^* : B \to C$ be an essentially
smooth $A$-algebra homomorphism, and let $(M, \rho) \in 
\msf{D}^{\mrm{b}}_{\mrm{f}}(\cat{Mod} B)_{\mrm{rig} / A}$.
\begin{enumerate}
\item The complex $f^{\sharp} M$ has an induced rigidifying 
isomorphism 
\[ f^{\sharp}(\rho) : f^{\sharp} M  \iso \opn{Sq}_{C / A}
f^{\sharp} M . \]
\item The rigid complex 
\[ f^{\sharp}(M, \rho) := 
\bigl( f^{\sharp} M, f^{\sharp}(\rho) \bigr) \in
\msf{D}^{\mrm{b}}_{\mrm{f}}(\cat{Mod} C)_{\mrm{rig} / A} \]
depends functorially on $(M, \rho)$ and on $f^*$. 
\end{enumerate}
\end{thm}

\begin{dfn} \label{dfn1.3}
Suppose $f^* : A \to B$ is essentially \'etale, so that
$f^{\sharp} M = B \otimes_A M$ for any 
$M \in \msf{D}(\cat{Mod} A)$. Let
$\mrm{q}^{\sharp}_{f; M} : M \to f^{\sharp} M$ be the morphism
$m \mapsto 1 \otimes m$. On the level of functors this gives a 
morphism
$\mrm{q}^{\sharp}_{f} : \bsym{1} \to f_* f^{\sharp}$
of functors from $\msf{D}(\cat{Mod} A)$ to itself.
\end{dfn}

In the situation of the definition, given 
$N \in \msf{D}(\cat{Mod} B)$, there is a canonical bijection
\[ \opn{Hom}_{\msf{D}(\cat{Mod} A)}(M, N) 
\iso \opn{Hom}_{\msf{D}(\cat{Mod} B)}(f^{\sharp} M, N) ,\
\phi \mapsto \bsym{1} \otimes \phi . \]
In particular, for 
$N := f^{\sharp} M$, the morphism $\mrm{q}^{\sharp}_{f; M}$
corresponds to the identity $\bsym{1}_{N}$.

\begin{dfn} \label{dfn1.4}
Let $A$ be a noetherian ring, 
let $A \to B$ be an essentially \'etale ring homomorphism, let
$M \in \msf{D}(\cat{Mod} A)$ and
$N \in \msf{D}(\cat{Mod} B)$.
A morphism $\phi : M \to N$ in $\msf{D}(\cat{Mod} A)$
is called a {\em nondegenerate 
localization morphism} if the corresponding morphism
$\bsym{1} \otimes \phi : f^{\sharp} M \to N$ in 
$\msf{D}(\cat{Mod} B)$ is an isomorphism. 
\end{dfn}

\begin{dfn} \label{dfn1.5}
Let $A$ be a noetherian ring, 
let $B$ and $C$ be essentially finite type $A$-algebras,
let $f^* : B \to C$ be an essentially \'etale 
$A$-algebra homomorphism, let
$(M, \rho) \in 
\msf{D}^{\mrm{b}}_{\mrm{f}}(\cat{Mod} B)_{\mrm{rig} / A}$
and let
$(N, \sigma) \in  
\msf{D}^{\mrm{b}}_{\mrm{f}}(\cat{Mod} C)_{\mrm{rig} / A}$.
A {\em rigid localization morphism} is a 
morphism $\phi : M \to N$ in 
$\msf{D}(\cat{Mod} B)$,
such that the corresponding morphism
$\bsym{1} \otimes \phi : 
f^{\sharp} (M, \rho) \to (N, \sigma)$ 
is a rigid morphism over $C$ relative to $A$.
\end{dfn}

\begin{prop} [{\cite[Proposition 6.8]{YZ5}}] \label{prop1.4}
Let $A$ be a noetherian ring, 
let $B$ and $C$ be essentially finite type $A$-algebras,
let $f^* : B \to C$ be an essentially \'etale 
$A$-algebra homomorphism, and
let $(M, \rho) \in 
\msf{D}^{\mrm{b}}_{\mrm{f}}(\cat{Mod} B)_{\mrm{rig} / A}$. 
Assume that $\mrm{RHom}_{B}(M, M) = B$.
Then the morphism $\mrm{q}^{\sharp}_{f; M}$
is the unique nondegenerate rigid localization morphism 
$(M, \rho) \to f^{\sharp}(M, \rho)$.
\end{prop}

The next result is about tensor products of rigid complexes.

\begin{thm}[{\cite[Theorem 0.4]{YZ5}}] \label{thm1.4}
Let $A$ be a noetherian ring, let $B, C$ be essentially finite 
type $A$-algebras, let $f^* : B \to C$ be an $A$-algebra
homomorphism, and let $(M, \rho) \in 
\msf{D}^{\mrm{b}}_{\mrm{f}}(\cat{Mod} B)_{\mrm{rig} / A}$
and
$(N, \sigma) \in 
\msf{D}^{\mrm{b}}_{\mrm{f}}(\cat{Mod} C)_{\mrm{rig} / B}$.
Assume that the canonical homomorphism 
$B \to \opn{Hom}_{\msf{D}(\cat{Mod} B)}(M, M)$
is bijective. Then the complex $M \otimes^{\mrm{L}}_{B} N$
is in $\msf{D}^{\mrm{b}}_{\mrm{f}}(\cat{Mod} C)$, it has finite 
flat dimension over $A$, and it has an 
induced rigidifying isomorphism 
\[ \rho \otimes \sigma : M \otimes^{\mrm{L}}_{B} N \iso
\opn{Sq}_{C / A} (M \otimes^{\mrm{L}}_{B} N) . \]
The rigid complex 
\[ (M, \rho) \otimes^{\mrm{L}}_{B} (N, \sigma) :=
(M \otimes^{\mrm{L}}_{B} N, \rho \otimes \sigma) \in
\msf{D}^{\mrm{b}}_{\mrm{f}}(\cat{Mod} C)_{\mrm{rig} / A} \]
depends functorially of $(M, \rho)$ and $(N, \sigma)$.
\end{thm}

Finally, a base change result.

\begin{thm}[{\cite[Theorem 6.9]{YZ5}}] \label{thm1.5}
Consider a commutative diagram of ring homomorphisms
\[ \UseTips \xymatrix @C=5ex @R=5ex {
A 
\ar[r]
& B 
\ar[r]^{f^*}
\ar[d]_{g^*}
& C 
\ar[d]^{h^*} \\
& B' 
\ar[r]_{f'^*}
& C' 
} \]
where $A$ is a noetherian ring, and 
$B, B', C$ and $C'$ are essentially finite type $A$-algebras.
Assume moreover that $g^* : B \to B'$ is a localization,
and the square is cartesian \tup{(}namely
$C' \cong B' \otimes_B C$\tup{)}. Let 
$(M, \rho) \in 
\msf{D}^{\mrm{b}}_{\mrm{f}}(\cat{Mod} B)_{\mrm{rig} / A}$,
let
$(N, \sigma) \in 
\msf{D}^{\mrm{b}}_{\mrm{f}}(\cat{Mod} C)_{\mrm{rig} / A}$,
and let 
\[ \phi : (N, \sigma) \to (M, \rho) \] 
be a rigid trace morphism over $B$ relative to $A$. 
Define 
$M' := g^{\sharp} M$ and $N' := h^{\sharp} N$.
There is a morphism
$\phi' : N' \to M'$ in $\msf{D}(\cat{Mod} B')$
gotten by composing the canonical isomorphism
$N' = C' \otimes_C N \cong B' \otimes_B N = g^{\sharp} N$
with
$g^{\sharp}(\phi) : g^{\sharp} N \to g^{\sharp} M = M'$.
So the diagram 
\[ \UseTips \xymatrix @C=5ex @R=5ex {
M
\ar[d]_{\mrm{q}^{\sharp}_{g; M}} 
& N
\ar[l]_(0.4){\phi} 
\ar[d]^{\mrm{q}^{\sharp}_{h; N}} \\
M' 
& N'
\ar[l]^(0.4){\phi'} 
} \]
is commutative. Then 
\[ \phi' : \bigl( N', h^{\sharp}(\sigma) \bigr) \to
\bigl( M', g^{\sharp}(\rho) \bigr)  \]
is a rigid trace morphism over $B'$ relative to $A$.

\end{thm}

Observe that there is no particular assumption on $f^*$.

\begin{rem}
The construction of the functor $\opn{Sq}_{B / A}$, and the proof 
of the theorems above in \cite{YZ5}, 
required heavy use of DG algebras. For the 
convenience of the reader we eliminated all reference to DG 
algebras in the definitions and statements in the present paper. 
However, we could not avoid using DG algebras in some of the 
proofs. 
\end{rem}

\section{Rigid Dualizing Complexes}
\label{sec.rigid-dual} 

In this section $\K$ is a fixed regular noetherian ring of finite 
Krull dimension. All algebras are by default essentially finite 
type $\K$-algebras, and all algebra homomorphisms are over $\K$.
We denote by $\cat{EFTAlg} / \K$ 
the category of essentially finite type $\K$-algebras.

Let us recall the definition of {\em dualizing complex} over a 
$\K$-algebra $A$ from \cite{RD}.
A complex $R \in \msf{D}^{\mrm{b}}_{\mrm{f}}(\cat{Mod} A)$
is called a dualizing complex if it has finite injective 
dimension, and the canonical morphism
$A \to \mrm{RHom}_{A}(R, R)$
in $\msf{D}(\cat{Mod} A)$ is an isomorphism.
It follows that the functor
$\mrm{RHom}_{A}(-, R)$
is an auto-duality (i.e.\ a contravariant equivalence) of 
$\msf{D}^{\mrm{b}}_{\mrm{f}}(\cat{Mod} A)$.
Note that since the ground ring $\K$ has 
finite global dimension, the complex $R$ has finite flat 
dimension over it. 

Following Van den Bergh \cite{VdB} we make the following 
definition.

\begin{dfn} \label{dfn2.1}
Let $A$ be an essentially finite type $\K$-algebra and let
$R$ be a dualizing complex over $A$. Suppose $R$ has a 
rigidifying isomorphism $\rho : R \iso \opn{Sq}_{A / \K} R$. 
Then the pair $(R, \rho)$ is 
called a {\em rigid dualizing complex over $A$ relative to $\K$}. 
\end{dfn}

By default all rigid dualizing complexes are relative to the 
ground ring $\K$. 

\begin{exa} \label{exa2.0}
Take the $\K$-algebra $A := \K$. 
The complex $R := \K$ is a dualizing complex 
over $\K$, since this ring is regular and finite dimensional. Let
\[ \rho^{\mrm{tau}} : \K \iso 
\mrm{RHom}_{\K \otimes_{\K} \K}(\K, \K \otimes_{\K} \K) =
\opn{Sq}_{\K / \K} \K \]
be the {\em tautological rigidifying isomorphism}. Then 
$(\K, \rho^{\mrm{tau}})$ is a 
rigid dualizing complex over $\K$ relative to $\K$.
\end{exa}

In \cite{VdB} it was proved that when $\K$ is a field,
a rigid dualizing complex $(R ,\rho)$ is unique up to isomorphism. 
And in \cite{YZ1} we proved that $(R, \rho)$ is in fact 
unique up to a unique rigid isomorphism (again, only when $\K$ is 
a field). These results are true in our setup too:

\begin{thm} \label{thm2.3} 
Let $\K$ be a regular finite dimensional noetherian ring,
let $A$ be an essentially finite type $\K$-algebra, and let 
$(R, \rho)$ be a rigid dualizing complex over $A$ 
relative to $\K$. Then $(R, \rho)$ is unique up to a 
unique rigid isomorphism.
\end{thm}

\begin{proof}
In view of \cite[Lemma 6.1]{YZ5} and \cite[Theorem 1.6]{YZ5}
we may assume that $\opn{Spec} A$ is connected. 
Suppose $(R', \rho')$ is another rigid dualizing complex over $A$. 
Then there is an isomorphism 
$R' \cong R \otimes_A L[n]$ for some invertible $A$-module $L$ and 
some integer $n$. Indeed 
$L[n] \cong \opn{RHom}_{A}(R, R')$; see \cite[Section V.3]{RD}.

Choose a K-flat DG algebra 
resolution $\K \to \til{A} \to A$ of $\K \to A$.
(If $\K$ is a field just take $\til{A} := A$.) So
\[ \begin{aligned}
& \opn{Sq}_{A / \K} R' \cong \opn{Sq}_{A / \K} (R_A \otimes_A L[n]) 
\\
& \quad = \mrm{RHom}_{\til{A} \otimes_{\K} \til{A}}
\bigl( A, (R_A \otimes_A L[n]) \otimes^{\mrm{L}}_{\K} 
(R_A \otimes_A L[n]) \bigr) \\
& \quad \cong^{\dag} \mrm{RHom}_{\til{A} \otimes_{\K} \til{A}}
( A, R_A \otimes^{\mrm{L}}_{\K} R_A )
\otimes^{\mrm{L}}_{A} L[n] \otimes^{\mrm{L}}_{A} L[n] \\
& \quad = (\opn{Sq}_{A / \K} R_A) 
\otimes^{\mrm{L}}_{A} L[n] \otimes^{\mrm{L}}_{A} L[n] 
\cong^{\diamondsuit} R_A \otimes_A L[n] \otimes_A L[n] . 
\end{aligned} \]
The isomorphism marked $\dag$ exists by 
\cite[Proposition 1.12]{YZ5} (with its condition (iii.b)),
and the isomorphism marked $\diamondsuit$ comes from 
$\rho : \opn{Sq}_{A / \K} R_A \iso R_A$. 
On the other and we have
$\rho' : R' \iso \opn{Sq}_{A / \K} R'$, which gives an isomorphism
\[ R_A \otimes_A L[n] \cong R_A \otimes_A L[n] \otimes_A L[n] . \]
Applying $\mrm{RHom}_A(R_A, -)$ to this isomorphism we get
$L[n] \cong L[n] \otimes_A L[n]$, and hence $L \cong A$ and 
$n = 0$. Thus we get an isomorphism $\phi_0 : R_A \iso R'$. 

The isomorphism $\phi_0$ might not be rigid, but there is some 
isomorphism $\phi_1$ making the diagram
\[ \UseTips \xymatrix @C=11ex @R=5ex {
R_A
\ar[r]^(0.4){\phi_1}
\ar[d]_(0.45){\rho_A}
&  R'
\ar[d]^(0.45){\rho'} \\
\opn{Sq}_{A / \K} R_A
\ar[r]^(0.5){\opn{Sq}_{\bsym{1}_A / \K}(\phi_0)}
& \opn{Sq}_{A / \K} R' } \]
commutative. Since 
$\mrm{Hom}_{\msf{D}(\cat{Mod} A)}(R_A, R')$
is a free $A$-module with basis $\phi_0$,
it follows that $\phi_1 = a \phi_0$ for some $a \in A^{\times}$. 
Then the isomorphism $\phi := a^{-1} \phi_0$ 
is the unique rigid isomorphism $R_A \iso R'$.
\end{proof}

In view of this result we are allowed to talk about {\em the} 
rigid dualizing complex over $A$ (if it exists).

Suppose $(M, \rho)$ is a rigid complex over $A$ relative to $\K$,
and $f^* : A \to B$ is a finite homomorphism of $\K$-algebras. 
Assume $f^{\flat} M$ has bounded cohomology. Then $f^{\flat} M$ 
has finite flat dimension over $\K$, and according to
Theorem \ref{thm1.2}(1) we get an induced rigid complex
$f^{\flat}(M, \rho)$ over $B$ relative to $\K$. 

\begin{prop}  \label{prop2.1}
Let $f^* : A \to B$ be a finite homomorphism of $\K$-algebras.
Assume a rigid dualizing 
complex $(R_A, \rho_A)$ over $A$ exists. Define
$R_B := f^{\flat}R_A \in \msf{D}(\cat{Mod} B)$
and $\rho_B := f^{\flat}(\rho_A)$. Then
$(R_B, \rho_B)$ is a rigid dualizing complex over $B$.
\end{prop}

\begin{proof}
The fact that $R_B$ is a dualizing complex over $B$ is proved in
\cite[Proposition V.2.4]{RD}. In particular $R_B$ has 
bounded cohomology.
\end{proof}

Suppose $(M, \rho)$ is a rigid complex over $A$ relative to $\K$,
and $f^* : A \to B$ is an essentially smooth homomorphism 
of $\K$-algebras. Then by Theorem \ref{thm1.3}(1) we get an 
induced rigid complex $f^{\sharp}(M, \rho)$ over $B$ relative 
to $\K$. 

\begin{prop}   \label{prop2.2}
Let $A$ be a $\K$-algebra, and 
assume $A$ has a rigid dualizing complex $(R_A, \rho_A)$.
Let $f^* : A \to B$ be an essentially smooth homomorphism. 
Define $R_B := f^{\sharp} R_A$
and $\rho_B := f^{\sharp}(\rho_A)$.
Then $(R_B, \rho_B)$ is a rigid dualizing complex over $B$.
\end{prop}

\begin{proof}
The complex $R_B$ is bounded. Hence to check it is dualizing is a 
local calculation on $\opn{Spec} B$. By 
\cite[Proposition 3.2(1)]{YZ5} 
we can assume that $A \to B$ is smooth. Now 
\cite[Theorem V.8.3]{RD} implies $R_B$ is dualizing.
\end{proof}

\begin{thm} \label{thm2.5}
Let $\K$ be a regular finite dimensional noetherian ring, and
let $A$ be an essentially finite type $\K$-algebra.
\begin{enumerate}
\item The algebra $A$ has a rigid dualizing complex 
$(R_A, \rho_A)$ relative to $\K$, 
which is unique up to a unique rigid isomorphism.
\item Given a finite homomorphism $f^* : A \to B$, there is a 
unique rigid isomorphism
\[ \phi^{\flat, \mrm{rig}}_f :
f^{\flat}(R_A, \rho_A) \iso (R_B, \rho_B) . \]
\item Given an essentially smooth homomorphism $f^* : A \to B$ , 
there is a unique rigid isomorphism
\[ \phi^{\sharp, \mrm{rig}}_f :
f^{\sharp}(R_A, \rho_A) \iso (R_B, \rho_B) . \]
\end{enumerate}
\end{thm}

\begin{proof}
(1) We can find algebras and homomorphisms
$\K \xar{f^*} C \xar{g^*} B \xar{h^*} A$,
where $C = \K[t_1, \ldots, t_n]$ is a polynomial algebra,
$g^*$ is surjective and $h^*$ is a localization. 
By Example \ref{exa2.0}, $(\K, \rho^{\mrm{tau}})$ 
is a rigid dualizing complex over $\K$. 
By Propositions \ref{prop2.1} and \ref{prop2.2} the complex
\[ R_A := h^{\sharp} g^{\flat} f^{\sharp} \K = A \otimes_B
\opn{RHom}_{C}(B, \Omega^n_{C / \K}[n]) \] 
is a rigid dualizing complex over $A$, with rigidifying 
isomorphism
$\rho_A := h^{\sharp} g^{\flat} f^{\sharp} (\rho^{\mrm{tau}})$.
Uniqueness was proved in Theorem \ref{thm2.3}.

\medskip \noindent
(2,3) Use Propositions \ref{prop2.1} and \ref{prop2.2} and Theorem 
\ref{thm2.3}.
\end{proof}

\begin{cor} \label{cor2.2}
Let $f^* : A \to B$ be a finite homomorphism 
in $\cat{EFTAlg} / \K$. 
There exists a unique nondegenerate rigid trace morphism
\[ \opn{Tr}_{f} = \opn{Tr}_{B / A}: 
(R_B, \rho_B) \to (R_A, \rho_A) . \]
\end{cor}

\begin{proof}
Since $f^{\flat} R_A$ is a dualizing complex over $B$ we know that 
\linebreak
$\opn{Hom}_{\msf{D}(\cat{Mod} B)}(f^{\flat} R_A, 
f^{\flat} R_A) = B$. 
So by \cite[Corollary 5.11]{YZ5}, 
$\opn{Tr}^{\flat}_{f; R_A} : f^{\flat} R_A \to R_A$
is the unique nondegenerate rigid
trace morphism between these two objects.

Composing $\opn{Tr}^{\flat}_{f; R_A}$
with the unique rigid isomorphism 
$R_B \cong f^{\flat} R_A$
guaranteed by Theorem \ref{thm2.3}, we get the unique rigid trace
$\opn{Tr}_{f} : R_B \to R_A$.
\end{proof}

\begin{cor} 
Let $A$ and $A'$ be in $\cat{EFTAlg} / \K$, 
with rigid dualizing complexes
$(R_A, \rho_A)$ and $(R_{A'}, \rho_{A'})$ respectively. 
Suppose $f^* : A \to A'$ is an essentially \'etale homomorphism.
Then there is exactly one nondegenerate rigid localization morphism 
\[ \mrm{q}_{f} = \mrm{q}_{A' / A} : (R_A, \rho_A) \to 
(R_{A'}, \rho_{A'}) . \]
\end{cor}

\begin{proof}
By Proposition \ref{prop2.2} we have a rigid complex
$f^{\sharp} R_A$ over $A'$, and by Proposition \ref{prop1.4}
there is a unique rigid localization morphism 
$\mrm{q}^{\sharp}_{f; R_A} : R_A \to f^{\sharp} R_A$.
According to Theorem \ref{thm2.3} there is a unique rigid 
isomorphism $f^{\sharp} R_A \cong R_{A'}$. By composing them we 
get the unique rigid localization map $\mrm{q}_{f}$.
\end{proof}

\begin{cor} \label{cor2.1}
Let $A, B$ and $A'$ be in $\cat{EFTAlg} / \K$, let
$f^* : A \to B$ be a finite homomorphism, and let 
$g^* : A \to A'$ be a localization. Define
$B' := A' \otimes_A B$, and let ${f'}^* : A' \to B'$ and 
$h^* : B \to B'$ be the induced homomorphisms. Then
\[ \opn{q}_{g} \circ\, \opn{Tr}_{f} = 
\opn{Tr}_{{f'}} \circ\, \opn{q}_{h} \]
in $\opn{Hom}_{\msf{D}(\cat{Mod} A)}(R_B, R_{A'})$.
\[ \UseTips \xymatrix @C=5ex @R=5ex {
A
\ar[r]^{f^*} 
\ar[d]_{g^*} 
& B
\ar[d]^{h^*} \\
A' 
\ar[r]_{{f'}^*}
& B'
} 
\qquad \quad
\UseTips \xymatrix @C=5ex @R=5ex {
R_A
\ar[d]_{\mrm{q}^{}_{g}} 
& R_B
\ar[l]_(0.4){\opn{Tr}_{f}} 
\ar[d]^{\mrm{q}^{}_{h}} \\
R_{A'}
& R_{B'}
\ar[l]^(0.4){\opn{Tr}_{f'}} 
} \]
\end{cor}

\begin{proof}
The $B'$-module 
$\opn{Hom}_{\msf{D}(\cat{Mod} A)}(R_B, R_{A'})$ is free of rank 
$1$, and both $\opn{q}_{g} \circ\, \opn{Tr}_{f}$ and
$\opn{Tr}_{{f'}} \circ\, \opn{q}_{h}$ are generators. 
So there is a unique invertible element $b' \in B'$ such that  
$\opn{Tr}_{{f'}} \circ\, \opn{q}_{h} = 
b' \opn{q}_{g} \circ\, \opn{Tr}_{f}$. 
Now by Theorem \ref{thm1.4} the morphism
\[ g^{\sharp}(\opn{Tr}_{f}) : h^{\sharp} (R_B, \rho_B) \to
g^{\sharp} (R_A, \rho_A)  \]
is a rigid trace morphism relative to $\K$. And so are
\[ \bsym{1} \otimes \opn{q}_g :
g^{\sharp} (R_A, \rho_A) \to (R_{A'}, \rho_{A'}) , \]
\[ \bsym{1} \otimes \opn{q}_h :
h^{\sharp} (R_B, \rho_B) \to (R_{B'}, \rho_{B'})  \]
and
\[ \opn{Tr}_{f'} :
(R_{B'}, \rho_{B'}) \to (R_{A'}, \rho_{A'}) 
. \]
We conclude that both
$(\bsym{1} \otimes \opn{q}_g) \circ g^{\sharp}(\opn{Tr}_{f})$
and
$\opn{Tr}_{f'} \circ (\bsym{1} \otimes \opn{q}_h)$
are nondegenerate rigid trace morphisms 
$h^{\sharp} (R_B, \rho_B) \to (R_{A'}, \rho_{A'})$
over $A'$ relative to $\K$, and therefore they must be equal
(cf.\ Theorem \ref{thm.rev.1}). But
$\opn{Hom}_{\msf{D}(\cat{Mod} A')}(h^{\sharp} R_B, R_{A'})$ 
is also a free $B'$-module of rank $1$, So $b' = 1$.
\end{proof}

Next comes a surprising result that basically
says ``all rigid complexes are dualizing''. 
The significance of this result is yet unknown. 

\begin{thm} \label{thm2.4}
Let $\K$ be a regular finite dimensional noetherian ring, and
let $A$ be an essentially finite type $\K$-algebra.
Assume $\opn{Spec} A$ is connected. Let $(M, \rho)$ be a nonzero 
rigid complex over $A$ relative to $\K$. Then $M$ is a dualizing 
complex over $A$. Hence there exists a unique rigid isomorphism 
$(M, \rho) \cong (R_A, \rho_A)$.
\end{thm}

The proof is after these two lemmas.

\begin{lem} \label{lem5.1}
Suppose $L \in \msf{D}^{\mrm{b}}_{\mrm{f}}(\cat{Mod} A)$
satisfies 
$\opn{Ext}^i_A(L, A / \mfrak{m}) = 0$
for all $i \neq 0$ and all maximal ideals $\mfrak{m} \subset A$.
Then $L$ is isomorphic to a finitely generated projective module
concentrated in degree $0$.
\end{lem}

\begin{proof}
This can be checked locally. Over a local ring $A_{\mfrak{m}}$, 
a minimal free resolution of $L_{\mfrak{m}}$ must have a single 
nonzero term, in degree $0$.
\end{proof}

\begin{lem} \label{lem5.2}
If $A$ is a field then the theorem is true, and moreover 
$M \cong A[d]$ for some integer $d$.
\end{lem}

\begin{proof}
Since any dualizing complex over the field $A$ is isomorphic to a 
shift of $A$, we see that the rigid dualizing complex 
satisfies $R_A \cong A[d]$ for some $d$. Therefore
$\opn{Sq}_{A / \K} (A[d]) \cong A[d]$. 

Consider the rigid complex $(M, \rho)$.
We can decompose 
$M \cong \boplus_{i = 1}^r A[d_i]$ in 
$\msf{D}(\cat{Mod} A)$. Then in the setup of Theorem 
\cite[Theorem 2.2]{YZ5} we have
\[ \begin{aligned}
\opn{Sq}_{A / \K} M 
& = \mrm{RHom}_{\til{A} \otimes_{\K} \til{A}}(A, 
M \otimes^{\mrm{L}}_{\K} M) \\
& \cong \boplus_{i,j}\
\mrm{RHom}_{\til{A} \otimes_{\K} \til{A}} \bigl( A, 
A[d_i] \otimes^{\mrm{L}}_{\K} A[d_j] \bigr) \\
& \cong \boplus_{i,j}\
\bigl( \opn{Sq}_{A / \K} (A[d]) \bigr) [d_i + d_j - 2d] . 
\end{aligned} \]
But $M \cong \opn{Sq}_{A / \K} M$, from which we obtain
\[ \boplus_{i = 1}^r A[d_i] 
\cong \boplus_{i,j = 1}^r \ A[d_i + d_j - d] . \]
Thus $r = 1$ and $d_1 = d$. 
\end{proof}

\begin{proof}[Proof of Theorem \tup{\ref{thm2.4}}]
Consider the object
$L := \opn{RHom}_A(M, R_A) \in 
\msf{D}^{\mrm{b}}_{\mrm{f}}(\cat{Mod} A)$.
Take an arbitrary maximal ideal $\mfrak{m} \subset A$, 
and define $B := A / \mfrak{m}$. We get a finite homomorphism 
$f^* : A \to B$. By Theorem \ref{thm1.2}(1) the complex
$\opn{RHom}_A(B, M) = f^{\flat} M$
is a rigid complex over $B$ relative to $\K$. So either 
$f^{\flat} M$ is zero, or, by Lemma \ref{lem5.2}, 
$f^{\flat} M \cong R_B \cong B[d]$ for some integer $d$.
On the other hand by Theorem \ref{thm2.5}(2) we have an 
isomorphism $f^{\flat} R_A \cong R_B$. 
Thus there are isomorphisms
\begin{equation} \label{eqn5.1}
\begin{aligned}
\opn{RHom}_A(L, B) & \cong \opn{RHom}_A(L, B[d])[-d] \\
& \cong \opn{RHom}_A \bigl( \opn{RHom}_A(M, R_A),  
\opn{RHom}_A(B, R_A) \bigr)[-d] \\
& \cong \opn{RHom}_A(B, M)[-d] 
= f^{\flat} M [-d] \cong B^r  
\end{aligned}
\end{equation}
where $r = 0,1$. In particular 
$\opn{Ext}^i_A(L, B) = 0$ for $i \neq 0$. 
Since $\mfrak{m}$ was arbitrary, Lemma \ref{lem5.1} tells us that 
we can assume $L$ is a finitely generated locally free $A$-module,
concentrated in degree $0$. Moreover, from the isomorphisms
(\ref{eqn5.1}) we see that the rank of $L$ at any point
of $\opn{Spec} A$ is at most $1$. Since $L$ is 
nonzero and $\opn{Spec} A$ is connected it follows that $L$ has 
constant rank $1$. 

At this stage we have $M \cong R_A \otimes_A L^{\vee}$,
where $L^{\vee} := \opn{Hom}_A(L, A)$.
We conclude that $M$ is a dualizing 
complex over $A$. By Theorem \ref{thm2.5}(1) there is a unique 
rigid isomorphism $R_A \cong M$.
\end{proof}

We end this section with some remarks and examples.

\begin{rem}
The assumption that the base ring $\K$ has finite global dimension 
seems superfluous. It is needed for technical reasons 
(bounded complexes have finite flat dimension), yet we don't 
know how to remove it. However, it seems necessary for $\K$ to be 
Gorenstein -- see next example. Also finiteness is important, as 
Example \ref{exa5.2} shows.
\end{rem}

\begin{rem}
A result of Van den Bergh (valid even in the noncommutative setup) 
says the following: if $A$ is a Gorenstein ring, then there is a 
canonical isomorphism 
\[ R_A \cong \opn{RHom}_A (\opn{Sq}_{A / \K} A, A ) \]
in $\msf{D}(\cat{Mod} A)$. 
\end{rem}

\begin{exa} \label{exa5.1}
Consider a field $\k$,
and let $\K = A := \k[t_1, t_2]/(t_1^2, t_2^2, t_1 t_2)$. 
Then $A$ does  not have a rigid dualizing complex relative to $\K$. 
The reason is that any dualizing complex over the artinian
local ring $A$ 
must be $R \cong A^*[n]$ for some integer $n$, where 
$A^* := \opn{Hom}_{\k}(A, \k)$. Now
$\opn{Sq}_{A / \K} R \cong R \otimes^{\mrm{L}}_{\K} R$, which has 
infinitely many nonzero cohomology modules. So there can be no 
isomorphism $R \cong \opn{Sq}_{A / \K} R$.
\end{exa}

\begin{exa} \label{exa5.2}
Take any field $\K$, and let $A := \K(t_1, t_2, \ldots)$,
the field of rational functions in countably many variables. So 
$A$ is a noetherian $\K$ algebra, but it is not of essentially 
finite type. Clearly $A$ has a dualizing complex (e.g.\ 
$R := A$), but as shown in \cite[Example 3.13]{YZ1}, there does not 
exist a rigid dualizing complex over $A$ relative to $\K$. 
\end{exa}

\begin{rem}
The paper \cite{SdS} by de Salas uses an idea similar to Van den 
Bergh's rigidity to define residues on local rings. However the 
results there are pretty limited. Lipman (unpublished notes) has 
an approach to duality using the fundamental class of the 
diagonal, which is close in spirit to the idea of rigidity.
\end{rem}

\section{The Twisted Inverse Image $2$-Functor}
\label{sec.2-func}  

In this section we translate properties of rigid 
dualizing complexes that were established in Section 
\ref{sec.rigid-dual} into properties of certain functors. 
As before we assume that $\K$ is
a regular noetherian ring of finite Krull dimension. 
All algebras are by default essentially finite type $\K$-algebras,
and all algebra homomorphisms are over $\K$.

Here is a review of the notion of $2$-functor,
following \cite[Section 5.15]{Ha}.
Let $\cat{Cat}$ be the $2$-category of all categories. 
The objects of $\cat{Cat}$ are the categories; the $1$-morphisms are 
the functors between categories; and the $2$-morphisms are the 
natural transformations between functors. Suppose 
$\cat{A}$ is some category. A $2$-functor (or pseudofunctor)
$F : \cat{A} \to \cat{Cat}$ is a triple
$F = (F_0, F_1, F_2)$ consisting of functions of the types 
explained below. The function $F_0$ is 
from the class of objects of $\cat{A}$ to 
the class of objects of $\cat{Cat}$; i.e.\ $F_0(A)$ is a category 
for any $A \in \cat{A}$. The function $F_1$ assigns to any morphism
$\alpha_0 \in \opn{Hom}_{\cat{A}}(A_0, A_1)$ a functor
$F_1(\alpha_0) : F_0(A_0) \to F_0(A_1)$. 
The function $F_2$ assigns to a composable pair of morphisms
$\alpha_0 \in \opn{Hom}_{\cat{A}}(A_0, A_1)$ and
$\alpha_1 \in \opn{Hom}_{\cat{A}}(A_1, A_2)$
a natural isomorphism
\[ F_2(\alpha_0, \alpha_1) :  F_1(\alpha_1 \circ \alpha_0) \iso
F_1(\alpha_1) \circ F_1(\alpha_0) \]
between functors
$F_0(A_0) \to F_0(A_2)$.
The data $(F_0, F_1, F_2)$ 
have to satisfy the compatibility condition
\[ F_2(\alpha_0, \alpha_1) \circ
F_2(\alpha_1 \circ \alpha_0, \alpha_2) = 
F_2(\alpha_1 , \alpha_2) \circ
F_2(\alpha_0, \alpha_2 \circ \alpha_1) \]
for any composable triple 
$A_0 \xar{\alpha_0} A_1 \xar{\alpha_1} A_2 \xar{\alpha_2} A_3$
of morphisms in $\cat{A}$.
Moreover for any object $A \in \cat{A}$, with identity morphism 
$\bsym{1}_A$, it is required that
$F_1(\bsym{1}_A) = \bsym{1}_{F_0(A)}$, the identity functor of the 
category $F_0(A)$; and also
$F_2(\bsym{1}_A, \bsym{1}_A)$ has to be 
the identity automorphism of the 
functor of $\bsym{1}_{F_0(A)}$.

We are going to construct a $2$-functor
$F : \cat{EFTAlg} / \K \to \cat{Cat}$.
Its $0$-component $F_0$ will assign the category
$F_0(A) := \msf{D}^{+}_{\mrm{f}}(\cat{Mod} A)$ to any algebra $A$.
The $1$-component $F_1$ will assign a functor
$F_1(f^*) = f^! : \msf{D}^{+}_{\mrm{f}}(\cat{Mod} A)
\to \msf{D}^{+}_{\mrm{f}}(\cat{Mod} B)$
to every algebra homomorphism $f^* : A \to B$. 
And for every composable pair of homomorphisms
$A \xar{f^*} B \xar{g^*} C$ there will be a natural isomorphism
$F_2(f^*, g^*) = \phi_{f, g}$.
Furthermore there will be isomorphisms of $2$-functors 
$f^! \cong f^{\flat}$ and $f^! \cong f^{\sharp}$ on suitable 
subcategories. These constructions will require a lot of notation,
and we tried to make this notation sensible. Isomorphisms
between functors of the same family (i.e.\ belonging to the same 
$2$-functor) will be labelled by ``$\phi$'' with modifiers
(e.g.\ $\phi_{f, g} : (f \circ g)^! \iso g^! f^!$).
Isomorphisms between functors belonging to different families will 
be labelled by ``$\psi$'' with modifiers 
(e.g.\ $\psi_f^{\sharp} : f^{\sharp} \iso f^!$).
When applying an isomorphism such as $\psi_{f}^{\sharp}$
to a particular object, say $M$, the 
notation will be 
$\psi_{f; M}^{\sharp} : f^{\sharp} M \iso f^! M$.

By Theorem \ref{thm2.5} any $\K$-algebra $A$ has a rigid 
dualizing complex $(R_A, \rho_A)$, which is unique up to a unique 
rigid isomorphism. For the sake of legibility we will often keep 
the rigidifying isomorphism $\rho_A$ implicit, and refer to the 
rigid dualizing complex $R_A$. 

\begin{dfn}
Given a $\K$-algebra $A$, with rigid dualizing complex $R_A$, 
define the {\em auto-duality functor} of 
$\msf{D}_{\mrm{f}}(\cat{Mod} A)$ relative to $\K$ to be
$\mrm{D}_A := \mrm{RHom}_{A}(-, R_A)$. 
\end{dfn}

Note that the functor $\mrm{D}_A$ exchanges the subcategories
$\msf{D}^{+}_{\mrm{f}}(\cat{Mod} A)$
and \linebreak
$\msf{D}^{-}_{\mrm{f}}(\cat{Mod} A)$.
Given a homomorphism of algebras $f^* : A \to B$ the functor
$\mrm{L} f^* = B \otimes^{\mrm{L}}_{A} -$ sends 
$\msf{D}^{-}_{\mrm{f}}(\cat{Mod} A)$ into
$\msf{D}^{-}_{\mrm{f}}(\cat{Mod} B)$.
This permits the next definition. 

\begin{dfn} \label{dfn4.5}
Let $f^* : A \to B$ be a homomorphism in $\cat{EFTAlg} / \K$. 
We define the {\em twisted inverse image} functor 
\[ f^! : \msf{D}^+_{\mrm{f}}(\cat{Mod} A) \to 
\msf{D}^+_{\mrm{f}}(\cat{Mod} B) \]
relative to $\K$ as follows. 
\begin{enumerate}
\rmitem{i} If $A = B$ and $f^* = \bsym{1}_A$ (the identity 
automorphism) then we let
$f^! := \bsym{1}_{\msf{D}^+_{\mrm{f}}(\cat{Mod} A)}$
(the identity functor). 
\rmitem{ii} Otherwise we define 
$f^! := \mrm{D}_B \, \mrm{L} f^* \, \mrm{D}_A$.
\end{enumerate}
\end{dfn}

Recall that for composable homomorphisms 
$A \xar{f^*} B \xar{g^*} C$ we sometimes write
$(f \circ g)^*$ instead of $g^* \circ f^*$. 

\begin{dfn} \label{dfn5.1}
Given two homomorphisms $A \xar{f^*} B \xar{g^*} C$
in $\cat{EFTAlg} / \K$ we define an isomorphism 
\[ \phi_{f, g} : (f \circ g)^! \iso g^! f^! \]
of functors 
$\msf{D}^{+}_{\mrm{f}}(\cat{Mod} A) 
\to \msf{D}^{+}_{\mrm{f}}(\cat{Mod} C)$
as follows. 
\begin{enumerate}
\rmitem{i} If either $A = B$ and $f^* = \bsym{1}_A$, or
$B = C$ and $g^* = \bsym{1}_B$, then $\phi_{f, g}$ is just the 
identity automorphism of
$(f \circ g)^! = g^! f^!$.
\rmitem{ii} Otherwise we use the adjunction isomorphism
$\bsym{1}_{\msf{D}^+_{\mrm{f}}(\cat{Mod} B)} \iso
\mrm{D}_B\, \mrm{D}_B$, 
together with the obvious isomorphism
$\mrm{L} (f \circ g)^* = \mrm{L} (g^* \circ f^*)
\cong \mrm{L} g^*\, \mrm{L} f^*$,
to obtain an isomorphism
\[ \begin{aligned}
(f \circ g)^! & = 
\mrm{D}_C\, \mrm{L} (f \circ g)^*\, \mrm{D}_A \cong
\mrm{D}_C\, \mrm{L} g^*\, \mrm{L} f^*\, \mrm{D}_A \\
&  \cong
\mrm{D}_C\, \mrm{L} g^*\, \mrm{D}_B\, \mrm{D}_B\,
\mrm{L} f^*\, \mrm{D}_A = g^! f^! .
\end{aligned} \]
\end{enumerate}
\end{dfn}

\begin{prop} \label{prop3.1}
For three homomorphisms 
$A \xar{f^*} B \xar{g^*} C \xar{h^*} D$ in 
$\cat{EFTAlg} / \K$ the isomorphisms $\phi_{-,-}$ satisfy
the compatibility condition 
\[ \phi_{g, h} \circ \phi_{f, g \circ h}  
 = \phi_{f, g} \circ \phi_{f \circ g, h} 
 : (f \circ g \circ h)^! \iso h^! g^! f^! . \]
Thus the assignment $f^* \mapsto f^!$ is the $1$-component 
of a $2$-functor $\cat{EFTAlg} / \K \to \cat{Cat}$, 
whose $0$-component is 
$A \mapsto \msf{D}^{+}_{\mrm{f}}(\cat{Mod} A)$. 
\end{prop}

In stating the proposition we were a bit sloppy with notation; for 
instance we wrote $\phi_{f, g}$, whereas the 
correct expression is 
$h^!(\phi_{f, g})$. This was 
done for the sake of legibility, and we presume the reader can 
fill in the omissions (also in what follows). 

\begin{proof}
By definition 
\[ (f \circ g \circ h)^! M = 
\mrm{D}_D \, \mrm{L} (f \circ g \circ h)^* \,
\mrm{D}_A M\,  \]
and
\[ h^! g^! f^! M = 
\mrm{D}_D\, \mrm{L} h^*\, \mrm{D}_C\, \mrm{D}_C\,
\mrm{L} g^*\, \mrm{D}_B\, \mrm{D}_B\,
\mrm{L} f^*\, \mrm{D}_A M . \] 
The two isomorphism 
$\phi_{g, h} \circ \phi_{f, g \circ h}$ and 
$\phi_{f, g} \circ \phi_{f \circ g, h}$
differ only by the order in which the adjunction isomorphisms
$\bsym{1}_{\msf{D}^+_{\mrm{f}}(\cat{Mod} B)} \cong
\mrm{D}_B\, \mrm{D}_B$
and
$\bsym{1}_{\msf{D}^+_{\mrm{f}}(\cat{Mod} C)} \cong
\mrm{D}_C\, \mrm{D}_C$ 
are applied, and correspondingly an isomorphism
$C \cong C \otimes^{\mrm{L}}_{B} B$ is replaced by
$D \cong D \otimes^{\mrm{L}}_{B} B$.
Due to standard identities the net effect is that 
$\phi_{g, h} \circ \phi_{f, g \circ h} =
\phi_{f, g} \circ \phi_{f \circ g, h}$.
\end{proof}

Suppose $A \xar{f^*} B \xar{g^*} C$ are finite 
homomorphisms in $\cat{EFTAlg} / \K$. 
Adjunction gives rise to an isomorphism
\[ g^{\flat} f^{\flat} M = 
\mrm{RHom}_{B} \bigl(C, \mrm{RHom}_{A}(B, M) \bigr) \cong
\mrm{RHom}_{A}(C, M) = (f \circ g)^{\flat} M \]
for any $M \in \msf{D}(\cat{Mod} A)$, and thus there is an 
isomorphism of functors
\[ \phi^{\flat}_{f, g} : (f \circ g)^{\flat} \iso
g^{\flat} f^{\flat} . \]
So $f^* \mapsto f^{\flat}$ is a $2$-functor on the subcategory of 
$\cat{EFTAlg} / \K$ consisting of all algebras, but only
finite homomorphisms.

\begin{thm} \label{thm3.2}
\begin{enumerate}
\item Let $f^* : A \to B$ be a finite homomorphism
in $\cat{EFTAlg} / \K$. The isomorphism 
$\phi^{\flat, \mrm{rig}}_f : f^{\flat} R_A \iso R_B$ of 
Theorem \tup{\ref{thm2.5}(2)} induces an isomorphism 
\[ \psi_f^{\flat} : f^{\flat} \iso f^! \]
of functors 
$\msf{D}^{+}_{\mrm{f}}(\cat{Mod} A) 
\to \msf{D}^{+}_{\mrm{f}}(\cat{Mod} B)$.
\item 
Given two finite homomorphism 
$A \xar{f^*} B \xar{g^*} C$, there is equality
\[ \psi_{g}^{\flat} \circ \psi_{f}^{\flat} \circ
\phi^{\flat}_{f,g} =
\phi_{f,g} \circ \psi_{f \circ g}^{\flat}  \]
as isomorphisms of functors
$(f \circ g)^{\flat} \iso g^{!} f^{!}$.
Thus the isomorphisms $\psi_{f}^{\flat}$ are $2$-functorial.  
\end{enumerate}
\end{thm}

\begin{proof}
(1) Take $M \in \msf{D}^{\mrm{b}}_{\mrm{f}}(\cat{Mod} A)$.
We then have a series of isomorphisms
\[ \begin{aligned} 
f^{!} M 
& = \mrm{RHom}_{B} \bigl( B \otimes^{\mrm{L}}_{A}
\mrm{RHom}_A(M, R_A), R_B \bigr) \\
& \cong^{\dag} 
\mrm{RHom}_{A} \bigl( \mrm{RHom}_A(M, R_A), R_B \bigr) \\
& \cong^{\diamondsuit} 
\mrm{RHom}_{A} \bigl( \mrm{RHom}_A(M, R_A),
\mrm{RHom}_A(B, R_A) \bigr) \\
& \cong^{\triangledown}
\mrm{RHom}_A(B, M) = 
f^{\flat} M ,
\end{aligned} \]
where the isomorphism marked $\dag$ is by the Hom-tensor 
adjunction; 
the isomorphism marked $\diamondsuit$ is induced by
$\phi^{\flat, \mrm{rig}}_f$; the isomorphism $\triangledown$ 
is because $R_A$ is a dualizing complex. 

\medskip \noindent(2) 
Given a second finite homomorphism $g^* : B \to C$ we obtain a 
commutative diagram of rigid isomorphisms
\[ \UseTips \xymatrix @C=12ex @R=5ex {
(f \circ g)^{\flat} R_A
\ar[d]_{\phi^{\flat}_{f,g; R_A}}
\ar[r]^{
\phi^{\flat, \mrm{rig}}_{f \circ g}}
&  R_C \\
g^{\flat}  f^{\flat} R_A
\ar[ur]_{\phi^{\flat, \mrm{rig}}_g \circ \phi^{\flat, \mrm{rig}}_f}
} \]
This is due to the fact that $\phi^{\flat}_{f,g; R_A}$ is a rigid 
isomorphism \cite[Theorem 5.3(2)]{YZ5}, together with
the uniqueness in Theorem \ref{thm2.5}(2).
From this, using standard identities, we deduce the desired 
equality.
\end{proof}

Next let $A \xar{f^*} B \xar{g^*} C$ be essentially smooth 
homomorphisms. Then according to \cite[Proposition 3.4]{YZ5}
there is an isomorphism of functors
$\phi^{\sharp}_{f,g} : (f \circ g)^{\sharp} \iso g^{\sharp}  
f^{\sharp}$.
The isomorphisms $\phi^{\sharp}_{-,-}$ are 
$2$-functorial on the subcategory of $\cat{EFTAlg} / \K$ 
consisting of essentially smooth homomorphisms. 

\begin{thm} \label{thm3.3}
\begin{enumerate}
\item Let $f^* : A \to B$ be an essentially smooth homomorphism in 
$\cat{EFTAlg} / \K$. The isomorphism 
$\phi^{\sharp, \mrm{rig}}_f : f^{\sharp} R_A \iso R_B$ of 
Theorem \tup{\ref{thm2.5}(3)} induces an isomorphism 
\[ \psi_f^{\sharp} : f^{\sharp} \iso f^! \]
of functors 
$\msf{D}^{+}_{\mrm{f}}(\cat{Mod} A) 
\to \msf{D}^{+}_{\mrm{f}}(\cat{Mod} B)$.
\item Given two essentially smooth homomorphisms
$A \xar{f^*} B \xar{g^*} C$, there is equality
\[ \psi_{g}^{\sharp} \circ \psi_{f}^{\sharp} \circ
\phi^{\sharp}_{f,g} =
\phi_{f,g} \circ \psi_{f \circ g}^{\sharp} \]
as isomorphisms of functors
$(f \circ g)^{\sharp} \iso g^{!} f^{!}$.
Thus the isomorphisms $\psi_{f}^{\sharp}$ 
are $2$-functorial.
\end{enumerate}
\end{thm}

\begin{proof}
(1) We can assume that $f^*$ has relative dimension $n$ (see 
\cite[Lemma 6.1]{YZ5}). 
Take any $M \in \msf{D}^{+}_{\mrm{f}}(\cat{Mod} A)$.
Then
\[ \begin{aligned} 
f^{!} M 
& = \mrm{RHom}_{B} \bigl( B \otimes^{\mrm{L}}_{A}
\mrm{RHom}_A(M, R_A), R_B \bigr) \\
& \cong^{\dag} 
\mrm{RHom}_{A} \bigl( \mrm{RHom}_A(M, R_A), R_B \bigr) \\
& \cong^{\diamondsuit} 
\mrm{RHom}_{A} \bigl( \mrm{RHom}_A(M, R_A),
\Omega^n_{B / A}[n] \otimes_A R_A \bigr) \\
& \cong^{\triangledown}
\mrm{RHom}_{A} \bigl( \mrm{RHom}_A(M, R_A),
  R_A \bigr) \otimes_A \Omega^n_{B / A}[n] \\
& \cong^{\ddag} 
M \otimes_A \Omega^n_{B / A}[n] = 
f^{\sharp} M ,
\end{aligned} \]
where the isomorphism marked $\dag$ is by the Hom-tensor 
adjunction; the isomorphism 
marked $\diamondsuit$ is induced by
$\phi^{\sharp, \mrm{rig}}_f$; the isomorphism $\triangledown$ 
is due to \cite[Proposition 1.12]{YZ5};
and the isomorphism $\ddag$ is by the adjunction isomorphism
$M \cong \opn{D}_A \opn{D}_A M$. 
We let 
$\psi_{f; M}^{\sharp} : f^{\sharp} M \iso f^! M$ 
be the composed isomorphism.

\medskip \noindent(2) 
Given a second essentially smooth 
homomorphism $g^* : B \to C$ we obtain a 
commutative diagram of rigid isomorphisms
\[ \UseTips \xymatrix @C=12ex @R=5ex {
(f \circ g)^{\sharp} R_A
\ar[d]_{\phi^{\sharp}_{f,g; R_A}}
\ar[r]^{
\phi^{\sharp, \mrm{rig}}_{f \circ g}}
&  R_C \\
g^{\sharp} f^{\sharp} R_A
\ar[ur]_{\phi^{\sharp, \mrm{rig}}_g \circ \phi^{\sharp, \mrm{rig}}_f}
} \]
This is due to the fact that $\phi^{\sharp}_{f,g; R_A}$ is a rigid 
isomorphism \cite[Theorem 6.3(3)]{YZ5}, together with
the uniqueness in Theorem \ref{thm2.5}(3).
From this, using standard identities, we deduce the desired 
equality.
\end{proof}

The next result explains the dependence of the 
twisted inverse image $2$-functor
$f \mapsto f^!$ on the base ring $\K$. Assume $\mbb{L}$ 
is an essentially finite type $\K$-algebra that's regular (but 
maybe not essentially smooth over $\K$). 
Just like for $\K$, any essentially 
finite type $\mbb{L}$-algebra $A$ has a rigid 
dualizing complex relative to $\mbb{L}$, which we denote by
$(R_{A / \mbb{L}}, \rho_{A / \mbb{L}})$. 
For any homomorphism $f^* : A \to B$ of $\mbb{L}$-algebras
there is a corresponding twisted inverse image functor 
$f^{! / \mbb{L}} : \msf{D}^{+}_{\mrm{f}}(\cat{Mod} A) \to
\msf{D}^{+}_{\mrm{f}}(\cat{Mod} B)$,
constructed using $R_{A / \mbb{L}}$ and $R_{B / \mbb{L}}$. 
Let $(R_{\mbb{L}}, \rho_{\mbb{L}})$ 
be the rigid dualizing complex of $\mbb{L}$ relative to $\K$. 

\begin{prop} 
Let $A$ be an essentially finite type $\mbb{L}$-algebra. Then 
$R_{\mbb{L}} \otimes^{\mrm{L}}_{\mbb{L}} R_{A / \mbb{L}}$
is a dualizing complex over $A$, and it has an induced 
rigidifying isomorphism relative to $\K$. Hence 
there is a unique isomorphism 
$R_{\mbb{L}} \otimes^{\mrm{L}}_{\mbb{L}} R_{A / \mbb{L}}
\cong R_A$
in $\msf{D}^{\mrm{b}}_{\mrm{f}}(\cat{Mod} A)_{\mrm{rig} / \K}$.
\end{prop}

\begin{proof}
We might as well assume $\opn{Spec} \mbb{L}$ is connected
(cf.\ \cite[Lemma 6.1]{YZ5}).
Since $\mbb{L}$ is regular, one has
$R_{\mbb{L}} \cong P[n]$
for some invertible $\mbb{L}$-module $P$ and some integer $n$.
Therefore 
$R_{\mbb{L}} \otimes^{\mrm{L}}_{\mbb{L}} R_{A / \mbb{L}}$ 
is a dualizing 
complex over $A$. According to Theorem \ref{thm1.4} the complex
$R_{\mbb{L}} \otimes^{\mrm{L}}_{\mbb{L}} R_{A / \mbb{L}}$ 
has an induced rigidifying 
isomorphism $\rho_{\mbb{L}} \otimes \rho_{A / \mbb{L}}$. 
Now use Theorem \ref{thm2.3}.
\end{proof}

\begin{cor} \label{cor3.1}
There is a canonical isomorphism
\[ (f^* \mapsto f^!) \cong (f^* \mapsto f^{! / \mbb{L}}) \]
 of $2$-functors $\cat{EFTAlg} / \mbb{L} \to \cat{Cat}$.
\end{cor}

\begin{proof}
The twist $R_{\mbb{L}} \otimes^{\mrm{L}}_{\mbb{L}} -$
gets canceled out in 
\[ f^{!} M \cong
\mrm{RHom}_{A} \bigl( \mrm{RHom}_A(M, 
R_{\mbb{L}} \otimes^{\mrm{L}}_{\mbb{L}}R_{A / \mbb{L}}), 
R_{\mbb{L}} \otimes^{\mrm{L}}_{\mbb{L}} R_{A / \mbb{L}} \bigr) . \]
\end{proof}

\begin{exa}
Take $\K := \mbb{Z}$ and $\mbb{L} := \mbb{F}_p = \mbb{Z} / (p)$
for some prime number $p$. Then 
$R_{\mbb{L}} = \mbb{L}[-1]$, and for any 
$A \in \cat{EFTAlg} / \mbb{L}$
we have $R_{A / \mbb{L}} \cong R_A[1]$.
\end{exa}

The final result of this section
connects our constructions to those of \cite{RD}. 
We shall restrict attention to the category 
$\cat{FTAlg} / \K$ of finite type $\K$-algebras.
Given a homomorphism
$f^* : A \to B$ let us denote by
\[ f^{! (\mrm{G})} : \msf{D}^{+}_{\mrm{f}}(\cat{Mod} A) \to
\msf{D}^{+}_{\mrm{f}}(\cat{Mod} B) \]
the twisted inverse image from \cite{RD}.

\begin{thm} \label{thm3.4}
Let $\K$ be a regular finite dimensional noetherian ring.
\begin{enumerate}
\item Given $A \in \cat{FTAlg} / \K$ let $\pi_A^* : \K \to A$ be the 
structural homomorphism, let 
$R_A^{(\mrm{G})} := \pi_A^{! (\mrm{G})} \K$,
and let $R_A$ be the rigid dualizing 
complex of $A$. Then there is an isomorphism
$R_A^{(\mrm{G})} \cong R_A$
in $\msf{D}(\cat{Mod} A)$.
\item There is an isomorphism 
$(f^* \mapsto f^!) \cong (f^* \mapsto f^{! (\mrm{G})})$
of $2$-functors \linebreak
$\cat{FTAlg} / \K \to \cat{Cat}$.
\end{enumerate}
\end{thm}

\begin{proof}
(1) Take homomorphisms $\K \xar{f^*} C \xar{g^*} B \xar{h^*} A$
as in the proof of Theorem \ref{thm2.5}(1).
Then 
$R_A \cong h^{\sharp} g^{\flat} f^{\sharp} \K$,
and also
$R_A^{(\mrm{G})} \cong h^{\sharp} g^{\flat} f^{\sharp} \K$.

\medskip \noindent
(2) For any $A$ fix an isomorphism 
$\psi^{(\mrm{G})}_A : R_A^{(\mrm{G})} \iso R_A$.
Let $\opn{D}^{(\mrm{G})}_A := \linebreak 
\opn{RHom}_A(-, R_A^{(\mrm{G})})$,
so we get an induced isomorphism
$\psi^{(\mrm{G})}_A : \opn{D}^{(\mrm{G})}_A \iso \opn{D}_A$
between the associated auto-duality functors. 
It is known that the $2$-functor $(f^* \mapsto f^{! (\mrm{G})})$
also satisfies 
$f^{! (\mrm{G})} \cong \opn{D}^{(\mrm{G})}_B \,
\mrm{L} f^* \, \opn{D}^{(\mrm{G})}_A$.
In this way we obtain an isomorphism of $2$-functors
$\psi^{(\mrm{G})} : f^{! (\mrm{G})} \iso f^!$. 
\end{proof}

\begin{rem} \label{rem3.11}
If $A$ is a flat $\K$-algebra, then flat base 
change for the theory in \cite{RD} 
endows $R^{(\mrm{G})}_A$ with a rigidifying 
isomorphism, thus making the isomorphism 
$R^{(\mrm{G})}_A \cong R_A$ canonical. (See \cite{Ye4}).
So in case $\K$ is a field one has a canonical isomorphism
$f^! \cong f^{! (\mrm{G})}$ of $2$-functors, and one may try to ask 
more precise questions, such as compatibility with the 
transformations in Theorems \ref{thm3.2} and \ref{thm3.3}.
\end{rem}

\section{Functorial Traces and Localizations}
\label{sec.funct-traces}

In this section again we work over a fixed base ring $\K$, which 
is assumed to be regular, noetherian and of finite Krull 
dimension. Recall that to each homomorphism 
$f^* : A \to B$ in $\cat{EFTAlg} / \K$ we constructed a
twisted inverse image functor
$f^! : \msf{D}^{+}_{\mrm{f}}(\cat{Mod} A) \to
\msf{D}^{+}_{\mrm{f}}(\cat{Mod} B)$. 

\begin{dfn} \label{dfn4.6}
Let $f^* : A \to B$ be a finite homomorphism in 
$\cat{EFTAlg} / \K$. 
Take any $M \in \msf{D}^{+}_{\mrm{f}}(\cat{Mod} A)$.
Then by Theorem \ref{thm3.2}(1) there is an isomorphism
$\psi_{f; M}^{\flat} : f^{\flat} M \iso f^{!} M$,
and by formula (\ref{eqn1.1}) there is a morphism
$\opn{Tr}^{\flat}_{f; M} : f^{\flat} M \to M$.
Define
\[ \opn{Tr}^{}_{f; M} := \opn{Tr}^{\flat}_{f; M} \circ\,
(\psi_{f; M}^{\flat})^{-1} : f^! M \to M . \]
On the level of functors this becomes a morphism 
\[ \opn{Tr}_{f} : f_* f^{!}  \to \bsym{1} \] 
of functors from $\msf{D}^{+}_{\mrm{f}}(\cat{Mod} A)$
to itself, called the {\em functorial trace map}.
\end{dfn}

\begin{prop} \label{prop.funct-traces.1}
\begin{enumerate}
\item The trace maps $\opn{Tr}^{}_{f}$ are nondegenerate. 
Namely, given a finite homomorphism $f^* : A \to B$ and an object
$M \in \msf{D}^{+}_{\mrm{f}}(\cat{Mod} A)$, the morphism 
$\opn{Tr}^{}_{f; M} : f^!M \to M$ is a nondegenerate trace 
morphism, in the sense of Definition \tup{\ref{dfn1.2}}.
\item The trace maps $\opn{Tr}^{}_{f}$ are $2$-functorial for 
finite homomorphisms. I.e.\ in the setup of Definition 
\tup{\ref{dfn5.1}}, 
with both $f^*, g^*$ finite, one has
\[ \opn{Tr}^{}_{f \circ g} = \opn{Tr}^{}_{f} \circ
\opn{Tr}^{}_{g} \circ\, \phi_{f,g} . \]
\end{enumerate}
\end{prop}

\begin{proof}
(1) This is true because 
$\opn{Tr}^{}_{f} : R_B \to R_A$ is a nondegenerate
trace morphism. See Corollary \ref{cor2.2}.

\medskip \noindent
(2) This follows from Theorem  \ref{thm3.2}(2).
\end{proof}

\begin{dfn} \label{dfn8.1}
Let $f^* : A \to B$ be an essentially \'etale homomorphism in 
\linebreak 
$\cat{EFTAlg} / \K$. Composing the localization map 
$\mrm{q}^{\sharp}_{f} : \bsym{1} \to f_* f^{\sharp}$
of Definition \ref{dfn1.3} with the isomorphism
$\psi_f^{\sharp} : f^{\sharp} \iso f^!$
of Theorem \ref{thm3.3}(1), we define the 
{\em functorial localization map}
\[ \mrm{q}^{}_{f} : \bsym{1} \to f_* f^{!} , \]
which is a morphism of functors from
$\msf{D}^{+}_{\mrm{f}}(\cat{Mod} A)$ to itself.
\end{dfn}

\begin{prop} \label{prop4.1}
\begin{enumerate}
\item The localization maps $\opn{q}^{}_{f}$ are nondegenerate. 
Namely, \linebreak given an essentially \'etale
homomorphism $f^* : A \to B$ and an object
$M \in \msf{D}^{+}_{\mrm{f}}(\cat{Mod} A)$, the morphism 
$\mrm{q}_{f;M} : M \to f^!M$ is a nondegenerate localization 
morphism, in the sense of Definition \tup{\ref{dfn1.4}}.
\item The localization maps $\opn{q}^{}_{f}$ are $2$-functorial for 
essentially \'etale homomorphisms. I.e.\ in the setup of Definition 
\tup{\ref{dfn5.1}}, with both $f^*, g^*$ essentially \'etale, 
one has
\[ \phi_{f,g} \circ \opn{q}^{}_{f \circ g} = 
\opn{q}^{}_{g} \circ \opn{q}^{}_{f} . \]
\end{enumerate}
\end{prop}

\begin{proof}
Assertion (1) is true because 
$\opn{q}^{\sharp}_{f; M} : M \to f^{\sharp} M$
is nondegenerate; see Proposition \ref{prop1.4}.
Assertion (2) is an immediate consequence 
of Theorem \ref{thm3.3}(2).
\end{proof}

\begin{prop}  \label{prop.funct-traces.3}
In the setup of Corollary \tup{\ref{cor2.1}} there is equality
\[ \mrm{q}_{g} \circ \opn{Tr}_{f} = 
\opn{Tr}_{f'} \circ\, \phi_{g, f'} \circ\, \phi_{f, h}^{-1}
 \circ\, \mrm{q}_{h} \]
of morphisms of functors
$f_* f^! \to g_* g^!$
from $\msf{D}^{+}_{\mrm{f}}(\cat{Mod} A)$ to itself.
\end{prop}

\begin{proof}
The functorial trace maps $\opn{Tr}_{f}$
are induced from the corresponding trace maps between the rigid 
dualizing complexes, via double dualization
(cf.\ Definition \ref{dfn4.6}). 
Likewise for the functorial
localization maps $\mrm{q}_{f}$. Thus the corollary is a consequence 
of Corollary \ref{cor2.1}.
\end{proof}

Here is an illustration of Proposition \ref{prop.funct-traces.3}.
Take $M \in \msf{D}^{+}_{\mrm{f}}(\cat{Mod} A)$,
and define
$N := f^! M$, $M' := g^! M$ and
$N' := h^! N$. Using the isomorphism 
$\phi_{g, f'} \circ\, \phi_{f, h}^{-1}$
we identify $N' = f'^! M'$. Then the diagram
\[ \UseTips \xymatrix @C=5ex @R=5ex {
M
\ar[d]_{\mrm{q}_{g}} 
& N
\ar[d]^{\mrm{q}_{h}}
\ar[l]_(0.4){\opn{Tr}_{f}} \\
M' 
& N' 
\ar[l]^(0.4){\opn{Tr}_{f'}} 
} \]
is commutative.

For a finite flat homomorphism $f^* : A \to B$ we 
denote by $\opn{tr}_{B/A} : B \to A$ the usual trace map.
Thus $\opn{tr}_{B/A}(b)$ is the trace of the operator $b$ acting 
by multiplication on the locally free $A$-module $B$.

\begin{thm} \label{thm3.1}
Suppose $f^* : A \to B$ is a finite \'etale homomorphism in
$\cat{EFTAlg} / \K$. Then for any 
$M \in \msf{D}^{\mrm{b}}_{\mrm{f}}(\cat{Mod} A)$
the diagram
\[ \UseTips \xymatrix @C=5ex @R=5ex {
B \otimes_A M
\ar[rr]^{\psi^{\sharp}_{f; M}} 
\ar[dr]_{\opn{tr}_{B/A} \otimes \bsym{1}_M \quad} 
& & f^! M 
\ar[dl]^{\opn{Tr}_{f; M}} \\
& M
} \]
is commutative.
\end{thm}

Note that 
\[ \psi^{\sharp}_{f; M}(b \otimes m) = 
b \cdot \mrm{q}_{f;M}(m) \]
for any $b \in B$ and $m \in M$. 

We need a few lemmas for the proof of the theorem. 

\begin{lem} \label{lem4.3}
Suppose $M, N, R, S \in \msf{D}^{\mrm{b}}_{\mrm{f}}(\cat{Mod} A)$,
with $R, S$ having finite injective dimension. Then there is an 
isomorphism
\[ M \otimes^{\mrm{L}}_{A} \mrm{RHom}_{A} \bigl(
\mrm{RHom}_{A}(N, R), S \bigr) \to
\mrm{RHom}_{A} \bigl(
\mrm{RHom}_{A}(M \otimes^{\mrm{L}}_{A} N, R), S \bigr) ,
\]
which is functorial in these four complexes, and, when $R = S$, 
it commutes with the adjunction morphisms
\[ M \otimes^{\mrm{L}}_{A} N \to 
M \otimes^{\mrm{L}}_{A} \mrm{RHom}_{A} \bigl(
\mrm{RHom}_{A}(N, R), R \bigr) \]
and
\[ M \otimes^{\mrm{L}}_{A} N \to 
\mrm{RHom}_{A} \bigl(
\mrm{RHom}_{A}(M \otimes^{\mrm{L}}_{A} N, R), R \bigr) . \]
\end{lem}

\begin{proof}
First consider $A$-modules $P, K, L$.
If $P$ is a finitely generated projective module, 
then the obvious homomorphisms of $A$-modules
\begin{equation} \label{eqn4.3}
P \otimes_{A} \mrm{Hom}_{A} (K, L) \to
\mrm{Hom}_{A} \bigl( \mrm{Hom}_{A}(P, A) \otimes_{A} K,
L \bigr)
\end{equation}
and
\begin{equation} \label{eqn4.4}
\mrm{Hom}_{A}(P, A) \otimes_{A} \mrm{Hom}_{A} (K, L) \to
\mrm{Hom}_{A} ( P \otimes_{A} K, L)
\end{equation}
are bijective. 

Now choose a resolution $P \to M$ by a
bounded above complex of finitely generated projective $A$-modules. 
Also choose resolutions $R \to I$ and $S \to J$ by 
bounded complexes of injective modules. 
So
\[ M \otimes^{\mrm{L}}_{A} \mrm{RHom}_{A} \bigl(
\mrm{RHom}_{A}(N, R), S \bigr) =
P \otimes_{A} \mrm{Hom}_{A} \bigl(
\mrm{Hom}_{A}(N, I), J \bigr) \]
and
\[ \mrm{RHom}_{A} \bigl(
\mrm{RHom}_{A}(M \otimes^{\mrm{L}}_{A} N, R), S \bigr) =
\mrm{Hom}_{A} \bigl(
\mrm{Hom}_{A}(P \otimes_{A} N, I), J \bigr) . \]
Because each $P^i$ is a finitely generated projective 
$A$-module, we have a bijection
\[ \begin{aligned}
& P^i \otimes_{A} \mrm{Hom}_{A} \bigl(
\mrm{Hom}_{A}(N^j, I^k), J^l \bigr) \\
& \qquad \cong 
\mrm{Hom}_{A} \bigl(
\mrm{Hom}_{A}(P^i, A) \otimes_{A} \mrm{Hom}_{A}(N^j, I^k) , 
J^l \bigr)
\end{aligned} \]
for any $i, j, k, l$. Indeed, this is exactly equation 
(\ref{eqn4.3}) with 
$K := \mrm{Hom}_{A}(N^j, I^k)$ and
$L := J^l$. Since the three 
complexes $N, I, J$ are bounded it follows that we have a canonical 
isomorphism of DG $A$-modules
\begin{equation}  \label{eqn4.5}
\begin{aligned}
& P \otimes_{A} \mrm{Hom}_{A} \bigl(
\mrm{Hom}_{A}(N, I), J \bigr) \\
& \qquad \cong 
\mrm{Hom}_{A} \bigl(
\mrm{Hom}_{A}(P, A) \otimes_{A} \mrm{Hom}_{A}(N, I) , 
J \bigr) .
\end{aligned} 
\end{equation}
Similarly, using equation (\ref{eqn4.4}), we have a canonical 
isomorphisms of DG $A$-modules
\begin{equation} \label{eqn4.6}
\begin{aligned}
& \mrm{Hom}_{A} \bigl(
\mrm{Hom}_{A}(P, A) \otimes_{A} \mrm{Hom}_{A}(N, I) , 
J \bigr) \\
& \qquad \cong 
\mrm{Hom}_{A} \bigl(
\mrm{Hom}_{A}(P \otimes_{A} N, I) , J \bigr) . \\
\end{aligned}
\end{equation}

Finally, is $S = R$ then we can take $J = I$, and then it is easy 
to track where $P \otimes_A N$ gets mapped in equations
(\ref{eqn4.5}) and (\ref{eqn4.6}).
\end{proof}

Given $M \in \msf{D}^{\mrm{b}}_{\mrm{f}}(\cat{Mod} A)$ let's write
\[ \chi_M := \opn{Tr}_{f; M} \circ\, \psi^{\sharp}_{f; M} :
f^{\sharp} M \to M \]
(this is a temporary definition). 
Since $f^{\sharp} A = B$ and $f^{\sharp} M = B \otimes_A M$,
there is a functorial isomorphism
$M \otimes^{\mrm{L}}_{A} f^{\sharp} A \iso f^{\sharp} M$.

\begin{lem} \label{lem4.1}
The morphism $\chi_M : f^{\sharp} M \to M$ is induced from 
$\chi_A : f^{\sharp} A \to A$ via the tensor operation
$M \otimes^{\mrm{L}}_{A} -$. Namely the diagram
\[ \UseTips \xymatrix @C=5ex @R=5ex {
M \otimes^{\mrm{L}}_{A} f^{\sharp} A
\ar[r]^(0.6){\cong}
\ar[d]_{\bsym{1}_M  \otimes \chi_A}
& f^{\sharp} M 
\ar[d]^{\chi_{M}} \\
M \otimes^{\mrm{L}}_{A} A
\ar[r]^(0.6){\cong}
& M 
} \]
is commutative. 
\end{lem}

\begin{proof}
Take any perfect complex
$N \in \msf{D}^{\mrm{b}}_{\mrm{f}}(\cat{Mod} A)$. 
Consider the sequence of morphisms
\begin{equation} \label{eqn4.1}
\begin{aligned}
& M \otimes^{\mrm{L}}_{A} f^{\sharp} N = 
M \otimes^{\mrm{L}}_A (B \otimes_A N) \\
& \quad \iso^{\dag} M \otimes^{\mrm{L}}_{A}
\Bigl( B \otimes_A \mrm{RHom}_{A} \bigl( \mrm{RHom}_A(N, R_A), 
R_A \bigr) \Bigr) \\
& \quad \iso^{\diamondsuit} M \otimes^{\mrm{L}}_{A}
\mrm{RHom}_{A} \bigl( \mrm{RHom}_A(N, R_A), 
R_B \bigr) \\
& \quad \to^{\triangledown}
 M \otimes^{\mrm{L}}_{A}
\mrm{RHom}_{A} \bigl( \mrm{RHom}_A(N, R_A), 
R_A \bigr) \\
& \quad \iso^{\dag} M \otimes^{\mrm{L}}_{A} N 
\end{aligned} 
\end{equation}
in which $\dag$ come from the adjunction isomorphism 
\[ N \cong \mrm{RHom}_{A} \bigl( \mrm{RHom}_A(N, R_A), 
R_A \bigr) ; \] 
$\diamondsuit$ comes from the localization map
$\mrm{q}_{f}  : R_A \to R_B$; and $\triangledown$ comes from the 
trace map $\opn{Tr}_{f} : R_B \to R_A$. 
By comparing these morphisms to the definition of 
$\psi^{\sharp}_{f; M}$ in Theorem \ref{thm3.3}(1), and the 
definition of $\opn{Tr}_{f; M}$ in Definition \ref{dfn4.6}, we see 
that the composition of all the morphisms in (\ref{eqn4.1})
is precisely
\[ \bsym{1}_M  \otimes \chi_N : M \otimes^{\mrm{L}}_{A} f^{\sharp} N
\to M \otimes^{\mrm{L}}_{A} N . \]

Now 
$M \otimes^{\mrm{L}}_{A}  N$ is also in
$\msf{D}^{\mrm{b}}_{\mrm{f}}(\cat{Mod} A)$,
because $N$ is perfect.
In parallel to the sequence of morphisms (\ref{eqn4.1}) there is 
another sequence
\begin{equation} \label{eqn4.2}
\begin{aligned}
& f^{\sharp}(M \otimes^{\mrm{L}}_{A}  N) = 
B \otimes_A M \otimes^{\mrm{L}}_A N \\
& \quad \iso^{\dag} 
B \otimes_A \mrm{RHom}_{A} 
\bigl( \mrm{RHom}_A(M \otimes^{\mrm{L}}_{A}N, R_A), R_A \bigr) \\
& \quad \iso^{\diamondsuit} 
\mrm{RHom}_{A} \bigl( \mrm{RHom}_A(M \otimes^{\mrm{L}}_{A} N, R_A), 
R_B \bigr) \\
& \quad \to^{\triangledown}
\mrm{RHom}_{A} \bigl( \mrm{RHom}_A(M \otimes^{\mrm{L}}_{A} N, R_A), 
R_A \bigr) \\
& \quad \iso^{\dag} M \otimes^{\mrm{L}}_{A} N .
\end{aligned} 
\end{equation}
The composition of all these morphisms is 
$\chi_{M \otimes^{\mrm{L}}_{A} N}$.

Since $A \to B$ is flat it follows that $R_B$ has finite injective 
dimension over $A$. 
According to Lemma \ref{lem4.3} at each step there is a 
canonical isomorphism from the object in
(\ref{eqn4.1}) to the corresponding object in
(\ref{eqn4.2}), and together these form a big commutative ladder. 
Therefore we get a commutative diagram
\[ \UseTips \xymatrix @C=5ex @R=5ex {
M \otimes^{\mrm{L}}_{A} f^{\sharp} N
\ar[r]^{\cong}
\ar[d]_{\bsym{1}_M  \otimes \chi_N}
& f^{\sharp} (M \otimes^{\mrm{L}}_{A} N) 
\ar[d]^{\chi_{M \otimes^{\mrm{L}}_{A} N}} \\
M \otimes^{\mrm{L}}_{A} N
\ar[r]^{=}
& M \otimes^{\mrm{L}}_{A} N
} \]
functorial in $M$ and $N$. Taking $N := A$ we get the desired 
assertion.
\end{proof}

\begin{lem} \label{lem4.2}
The morphism
\[ \chi_A : f^{\sharp} (A, \rho^{\mrm{tau}}_A) \to
(A, \rho^{\mrm{tau}}_A) \]
is a nondegenerate rigid trace morphism relative to $A$.
\end{lem}

\begin{proof}
Since $f^{\sharp} \cong f^{\flat}$ and
$\opn{Tr}_{f; A}^{\flat} : f^{\flat} A \to A$
is a nondegenerate rigid trace morphism relative to $A$, it 
follows that $b \chi_A$ is a nondegenerate rigid trace morphism
relative to $A$ for some unique $b \in B^{\times}$. We are going 
to prove that $b = 1$.

Take the rigid dualizing complex
$(R_A, \rho_A) \in 
\msf{D}^{\mrm{b}}_{\mrm{f}}(\cat{Mod} A)_{\mrm{rig} / \K}$.
Recall the tensor operation for rigid complexes (Theorem 
\ref{thm1.4}). Trivially the isomorphism 
\[ (R_A, \rho_A) \otimes^{\mrm{L}}_{A} 
(A, \rho^{\mrm{tau}}) \cong (R_A, \rho_A) \]
is rigid over $A$ relative to $\K$. According to 
\cite[Theorem 6.3(2)]{YZ5} the isomorphism 
\[ (R_A, \rho_A) \otimes^{\mrm{L}}_{A} 
f^{\sharp} (A, \rho^{\mrm{tau}}) \cong 
f^{\sharp} (R_A, \rho_A) \]
is rigid over $B$ relative to $\K$. 
Since the tensor operation of rigid 
complexes is functorial, it follows that the morphism 
\[ \bsym{1}_{R_A} \otimes b \chi_A : f^{\sharp} (R_A, \rho_A) \to
(R_A, \rho_A) \]
is a nondegenerate rigid trace morphism. 
Next, from Lemma \ref{lem4.1} we know that
\[ \bsym{1}_{R_A} \otimes \chi_A = \chi_{R_A} :
f^{\sharp} R_A \to R_A . \]
We conclude that 
\[ b \chi_{R_A} : f^{\sharp} (R_A, \rho_A) \to
(R_A, \rho_A) \]
nondegenerate rigid trace morphism over $A$ relative to $\K$.  

On the other hand
\[ \chi_{R_A} =  
\opn{Tr}_{f}^{\flat} \circ\, 
(\psi^{\flat, \mrm{rig}}_{f})^{-1} \circ 
\psi^{\sharp, \mrm{rig}}_{f} : f^{\sharp} (R_A, \rho_A) \to
(R_A, \rho_A) , \]
which itself is a rigid trace morphism over $A$ relative to $\K$. 
Since there is only one
nondegenerate rigid trace morphism 
$f^{\sharp} (R_A, \rho_A) \to (R_A, \rho_A)$, 
it follows that $b \chi_{R_A} = \chi_{R_A}$, and hence $b = 1$.
\end{proof}

\begin{proof}[Proof of Theorem \tup{\ref{thm3.1}}]
By definition $\opn{tr}_{B/A} \otimes\, \bsym{1}_M$
is induced from $\opn{tr}_{B/A}$. And according to Lemma \ref{lem4.1}
the morphism 
$\chi_M = \opn{Tr}_{f; M} \circ \, \psi^{\sharp}_{f; M}$
is also induced from 
$\chi_A$. Therefore it is enough to look at 
$M = A$. We must show that
$\opn{Tr}_{f; A} \circ\, \psi^{\sharp}_{f; A} = 
\opn{tr}_{B / A}$.

Let $B^{\mrm{e}} := B \otimes_A B$. 
According to \cite[Proposition 3.15]{YZ5} we know that there is a 
canonical ring isomorphism
$B^{\mrm{e}} \cong B \times B'$, where the factor $B'$ is the 
kernel of the multiplication map 
$B^{\mrm{e}} \to B$. Thus the surjective
$B^{\mrm{e}}$-module homomorphism
$B^{\mrm{e}} \to B$ has a canonical splitting 
$\nu : B \to B^{\mrm{e}}$. 
 From the proofs of \cite[Theorems 6.3(2) and 3.14(3)]{YZ5}
we see that the rigidifying isomorphism 
$\rho := f^{\sharp}(\rho^{\mrm{tau}})$ of
$B = f^{\sharp} A$ is precisely
\[ \nu : B \iso \opn{Hom}_{B^{\mrm{e}}}(B, B^{\mrm{e}})
= \opn{Sq}_{B / A} B . \]
Now Lemma \ref{lem4.2} says that the morphism
\[ \chi_A =  \opn{Tr}_{f; A} \circ\, \psi^{\sharp}_{f; A}
: (B, \rho) \to (A, \rho^{\mrm{tau}}) \]
is a nondegenerate rigid trace morphism relative to $A$. 
Since $\opn{tr}_{B / A}$ is also nondegenerate, 
it suffices to prove that 
$\opn{tr}_{B / A}$ is a rigid trace morphism relative to $A$. 

We have finite flat ring homomorphisms
$A \to B \xar{g^*} B^{\mrm{e}}$,
where $g^*$ is $b \mapsto b \otimes 1$. 
Because $B^{\mrm{e}} \cong B \times B'$ we have
\[ \opn{tr}_{B^{\mrm{e}} / B} \bigl( \nu(b) \bigr) = 
\opn{tr}_{(B \times B') / B}(b, 0) = b \] 
for any $b \in B$. But
$\opn{tr}_{B^{\mrm{e}} / A} = 
\opn{tr}_{B / A} \otimes \opn{tr}_{B / A}$,
and we know that $\rho = \nu$. Also the traces are transitive. So
\[ \opn{tr}_{B / A}(b) = 
(\opn{tr}_{B / A} \circ \opn{tr}_{B^{\mrm{e}} / B}) \bigl(
\nu(b) \bigr) = 
\opn{tr}_{B^{\mrm{e}} / A} \bigl( \nu(b) \bigr) = 
(\opn{tr}_{B / A} \otimes \opn{tr}_{B / A}) \bigl( \rho(b) \bigr) 
. \]
This means that indeed $\opn{tr}_{B / A}$ is a rigid trace
morphism relative to $A$.
\end{proof}

\section{Traces of Differential Forms}
\label{sec.traces}

A useful feature of Grothendieck duality
theory is that it gives rise to traces of differential forms. Such 
traces are quite hard to construct directly
(cf.\ \cite{Li1}, \cite{Hu} and \cite{Ku}). 
The aim of this section is to construct trace maps and to study 
some of their properties. The connection of our 
constructions to \cite{RD} is via Theorem \ref{thm3.4}.

As before $\K$ is a regular noetherian ring of finite Krull 
dimension. 

\begin{dfn} \label{dfn.traces.1}
Suppose $A \xar{f^*} B \xar{g^*} C$ are homomorphisms in 
$\cat{EFTAlg} / \K$, with $f^* : A \to B$ and 
$(f \circ g)^* : A \to C$ 
essentially smooth of relative dimension $n$, and 
$g^* : B \to C$ finite.
According to Theorem \ref{thm3.3}(1) there are isomorphisms
$\psi^{\sharp}_{f; A} : \Omega^n_{B / A}[n] \iso f^! A$
and
$\psi^{\sharp}_{f \circ g; A} : \Omega^n_{C / A}[n] \iso 
(f \circ g)^! A$.
 From Definition \ref{dfn4.6} there is a trace map
$\opn{Tr}_g : g_* g^! \iso \bsym{1}$,
and from Definition \ref{dfn5.1} there is an isomorphism
$\phi_{f, g} : (f \circ g)^! \iso g^! f^!$.
Define
\[ \opn{Tr}_{C / B / A}  = \opn{Tr}_{f / A} := 
(\psi^{\sharp}_{f})^{-1} \circ \opn{Tr}_g \circ\, \phi_{f, g} \circ
\psi^{\sharp}_{f \circ g} [-n] . \]
Thus
\[ \opn{Tr}_{C / B / A} : \Omega^n_{C/A} \to \Omega^n_{B/A} \]
is a $B$-linear homomorphism
called the {\em trace map}.
\end{dfn}

\begin{thm} \label{thm5.1}
Let $A \to B \to C$ be homomorphisms in 
$\cat{EFTAlg} / \K$ as in Definition \tup{\ref{dfn.traces.1}}.
\begin{enumerate}
\item The trace map
$\opn{Tr}_{C / B / A} : \Omega^n_{C/A} \to \Omega^n_{B/A}$
is nondegenerate, i.e.\ it induces a bijection
$\Omega^n_{C/A} \cong \opn{Hom}_{B}(C, \Omega^n_{B/A})$.
\item Suppose that $C \to D$ is a finite
homomorphism, such that the composed homomorphism
$A \to D$ is also essentially smooth 
of relative dimension $n$. Then there is equality 
\[  \opn{Tr}_{D / B / A} = \opn{Tr}_{C / B / A}  \circ
\opn{Tr}_{D / C / A} . \]
\end{enumerate}
\end{thm}

\begin{proof}
These assertions follow directly
from Proposition \ref{prop.funct-traces.1}(1,2)
respectively. 
\end{proof}

\begin{rem} \label{rem5.1}
In the setup of Definition \ref{dfn.traces.1}, 
the homomorphism $g^* : B \to C$ is actually 
flat. The proof will be published elsewhere; but here is 
the idea. It suffices to check the flatness of $g^*$ 
after passing to the induced homomorphisms
$\what{A}_{\mfrak{p}} \to \what{B}_{\mfrak{q}}
\to \what{C}_{\mfrak{r}}$ 
between the complete local rings, where 
$\mfrak{r} \in \opn{Spec} C$ is arbitrary, 
$\mfrak{q} := g(\mfrak{r})$ and $\mfrak{p} := f(\mfrak{q})$.
We then prove that these homomorphisms can 
be lifted to homomorphisms 
$\til{A} \to \til{B} \to \til{C}$, 
such that all three rings are regular
complete local rings, and $\til{B} \to \til{C}$ 
is finite and injective. It is now a classical result that 
$\til{B} \to \til{C}$ is flat. 
\end{rem}

\begin{rem}
R. H\"ubl has communicated to us that results of Kunz \cite{Ku} 
imply that the homomorphism $g^* : B \to C$ above is not only 
flat, but in fact a locally complete intersection. Hence according 
to \cite[Section 16]{Ku} there is a trace map
$\sigma_{C/B} : \Omega_{C / A} \to \Omega_{B / A}$,
which is a homomorphism of DG modules. In degree 
$n$ it is a nondegenerate $A$-linear map
$\sigma_{C/B}^n : \Omega^n_{C / A} \to \Omega^n_{B / A}$.
Presumably Kunz's trace map $\sigma_{C/B}^n$
coincides with our trace map $\opn{Tr}_{C / B / A}$, although this 
is quite hard to verify (cf.\ Propositions 
\ref{prop5.1} and \ref{prop5.2} below). 
\end{rem}

\begin{rem}
One can show that 
\[ \opn{Tr}_{C / B / A}[n] : (f \circ g)^{\sharp} 
(A, \rho^{\mrm{tau}}) \to
f^{\sharp} (A, \rho^{\mrm{tau}}) \]
is the unique nondegenerate rigid trace morphism over $B$
relative to $A$ between these two rigid complexes. This means that 
the trace map $\opn{Tr}_{C / B / A}$ is actually independent of 
base ring $\K$. Cf.\ a similar phenomenon for the 
cohomological residue map in \cite{Ye4}. 
\end{rem}

\begin{dfn} \label{dfn.traces.2}
Let $A \xar{f^*} B \xar{g^*} C$ be homomorphisms in
$\cat{EFTAlg} / \K$,
with $f^* : A \to B$ essentially smooth of relative dimension $n$, 
and $g^* : B \to C$ 
essentially \'etale. Define a $B$-linear homomorphism
\[ \opn{q}_{C / B / A} : \Omega^n_{B / A} \to \Omega^n_{C / A} \]
by the formula
\[ \opn{q}_{C / B / A} 
\bigl( b_0 \d b_1 \wedge \cdots \wedge \d b_n \bigr)
= g^*(b_0) \d \bigl( g^*(b_1) \bigr) \wedge \cdots \wedge 
\d \bigl( g^*(b_n) \bigr)  \]
for any $b_0, \ldots, b_n \in B$. Here $\d$ is the de Rham 
differential.  
\end{dfn}

It is trivial that $\opn{q}_{C / B / A}$ is a nondegenerate
localization homomorphism, namely
\[ \bsym{1} \otimes \opn{q}_{C / B / A} : 
C \otimes_B \Omega^n_{B / A} \to \Omega^n_{C / A} \]
is bijective.

\begin{lem} \label{lem5.3}
Let $A \xar{f^*} B \xar{g^*} C$ be as in Definition 
\tup{\ref{dfn.traces.2}}. Then
\[ \opn{q}_{C / B / A}[n] : f^{\sharp} (A, \rho^{\mrm{tau}}) 
\to (f \circ g)^{\sharp}(A, \rho^{\mrm{tau}})  \]
is the unique nondegenerate rigid localization 
morphism over $B$ relative to $A$ between these two rigid 
complexes.
\end{lem}

\begin{proof}
By \cite[Theorem 6.3(3)]{YZ5} the obvious isomorphism
\[ \nu : g^{\sharp} f^{\sharp} (A, \rho^{\mrm{tau}}) \iso
(f \circ g)^{\sharp} (A, \rho^{\mrm{tau}}) \]
is rigid. And according to \cite[Proposition 6.8]{YZ5}
the localization map
\[ \opn{q}^{\sharp}_g :  f^{\sharp} (A, \rho^{\mrm{tau}}) \to
g^{\sharp} f^{\sharp} (A, \rho^{\mrm{tau}}) \]
is the unique nondegenerate rigid localization map. 
But 
$\opn{q}_{C / B / A}[n] = \nu \circ \opn{q}^{\sharp}_g$.
\end{proof}

The next result says that the trace maps commute with 
localizations. 

\begin{prop} \label{prop5.3} 
Consider a commutative diagram
\[ \UseTips \xymatrix @C=5ex @R=5ex {
A 
\ar[r]
& B 
\ar[r]
\ar[d]
& C 
\ar[d] \\
& B' 
\ar[r]
& C' 
} \]
in $\cat{EFTAlg} / \K$, in which the homomorphisms $A \to B$ and
$A \to C$ are essentially smooth of relative 
dimension $n$; $B \to C$ is finite; $B \to B'$ is a localization; 
and the square is cartesian \tup{(}i.e.\ 
$C' \cong B' \otimes_B C$\tup{)}. Then the diagram
\[ \UseTips \xymatrix @C=11ex @R=5ex {
\Omega^n_{B / A}
\ar[d]_{\opn{q}_{B' / B / A}}
& \Omega^n_{C / A}
\ar[l]_(0.45){\opn{Tr}_{C / B / A}} 
\ar[d]^{\opn{q}_{C' / C / A}} \\
\Omega^n_{B' / A}
& \Omega^n_{C' / A}
\ar[l]^(0.45){\opn{Tr}_{C' / B' / A}}
} \]
is commutative.
\end{prop}

\begin{proof}
Let's write
$e^* : A \to B$ and $g^* : B \to B'$.
By Lemma \ref{lem5.3} the localization map
$\opn{q}_{B' / B / A}[n] : \Omega^n_{B / A}[n] 
\to \Omega^n_{B' / A}[n]$
gets sent to the localization map
$\opn{q}_g : e^! A \to (e \circ g)^! A$
under the isomorphisms 
$\psi^{\sharp}_{e; A} : \Omega^n_{B / A}[n]  \iso e^! A$
and
$\psi^{\sharp}_{e \circ g; A} : \Omega^n_{B / A}[n]  \iso 
(e \circ g)^! A$
of Theorem \ref{thm3.3}. Likewise for 
$\opn{q}_{C' / C / A}[n]$. Now we can use Proposition 
\ref{prop.funct-traces.3}. 
\end{proof}

To finish the paper here are two nice properties of the trace 
map.

\begin{prop} \label{prop5.1} 
In the setup of Definition \tup{\ref{dfn.traces.1}},  
assume that the ring $A$ is reduced. Then for any 
$\beta \in \Omega^n_{B / A}$ and $c \in C$ one has
\[ \opn{Tr}_{C/B/A}(c \beta) = \opn{tr}_{C/B}(c) \cdot \beta \in
\Omega^n_{B / A} . \]
\end{prop}

\begin{proof}
Denote by $A'$ the total ring of fractions of $A$, and let
$B' := A' \otimes_A B$. Since $A \to A'$ is injective, so is
$\Omega^n_{B / A} \to \Omega^n_{B' / A}$. 
Now $A'$ is a finite product of fields, and hence $B'$ is a finite 
product of integral domains (cf.\ 
\cite[Proposition 3.2]{YZ5}). 
Let $B''$ be the total ring of fractions of $B'$. Then $B''$ is a 
finite product of fields, and 
$\Omega^n_{B' / A} \to \Omega^n_{B'' / A}$
is injective. Note that $C' := A' \otimes_A C$
is also a finite product of integral domains, so 
$C'' := B'' \otimes_B C'$ is a finite product of fields.
Due to Proposition \ref{prop5.3}
we may replace $A \to B \to C$
with $A \to B'' \to C''$. We can also localize at one of the 
factors of $B''$. Thus we might as well assume that $B$ is a field 
and $C = \prod C_i$ is a finite product of fields.

Since each homomorphism $C \to C_i$ is both finite and a 
localization, Theorem \ref{thm5.1}(2) and Proposition \ref{prop5.3}
imply that for any 
$\gamma \in \Omega^n_{C_i / A}$ the $i$-th component of
$\opn{Tr}_{C_i / C / A}(\gamma) \in \Omega^n_{C / A} = 
\boplus_j \Omega^n_{C_j / A}$
is $\gamma$, and all other components are zero.
We conclude that we can replace $C$ with $C_i$; i.e.\ $C$ can be 
assumed to be a field.

Now there are two cases to look at. If the finite field extension
$B \to C$ is separable then 
$\opn{Tr}_{C/B/A}(c \beta) = \opn{tr}_{C/B}(c) \cdot \beta$
by Theorem \ref{thm3.1}.
 On the other hand, if $B \to C$ is 
inseparable then $\opn{tr}_{C/B}(c) = 0$ and also
$\mrm{q}_{C / B / A}(\beta) = 0$.
\end{proof}

\begin{prop} \label{prop5.2} 
Let $A \in \cat{EFTAlg} / \K$ be an integral domain whose field of 
fractions has characteristic $0$. Let
$B := A[s]$ and $C := A[t]$ be polynomial algebras in one variable 
each. Define an $A$-algebra homomorphism $f^* : B \to C$ by
$f^*(s) := t^n$ for some positive integer $n$. Then
\[ \opn{Tr}_{C / B / A} \bigl( t^{n-1} \d t \bigr) = \d s , \]
and
\[ \opn{Tr}_{C / B / A} \bigl( t^i \d t \bigr) = 0
\ \text{ for }\ 0 \leq i \leq n-2 . \]
\end{prop}

\begin{proof}
Let $A'$ be the fraction field of $A$. And let's write
$\opn{Tr}_{C / B / A} \bigl( t^i \d t \bigr) = p_i(s) \d s$
with $p_i(s) \in B = A[s]$.
So we need to prove that $p_{n-1}(s) = 1$ and $p_i(s) = 0$
for $0 \leq i \leq n-2$. 
Due to Proposition \ref{prop5.3}
we can localize to $B' := A'[s]$ and $C' := A'[t]$. 
Denote by ${f'}^* : B' \to C'$ the corresponding homomorphism.
Take any 
nonzero $\lambda \in A'$. There are $A'$-algebra automorphisms
$g^*_{\lambda} : B' \to B'$ and $h^*_{\lambda} : C' \to C'$,
defined by $g^*_{\lambda}(s) := \lambda^n s$ and
$h^*_{\lambda} (t) := \lambda t$,
and these satisfy 
$h^*_{\lambda} \circ {f'}^* = {f'}^* \circ g^*_{\lambda}$.
Since the trace is functorial 
(Theorem \ref{thm5.1}(2)) we have
$\opn{Tr}_{f' / A} \circ \opn{Tr}_{h_{\lambda} / A} =
\opn{Tr}_{g_{\lambda} / A} \circ \opn{Tr}_{f' / A}$.
But by Proposition \ref{prop5.1},
$\opn{Tr}_{g_{\lambda} / A} 
\bigl( g_{\lambda}^*(\beta) \bigr) = \beta$
for any $\beta \in \Omega^n_{B' / A}$; so that
$\opn{Tr}_{g_{\lambda} / A} \bigl( p(s) \d s \bigr) = 
p(\lambda^{-n} s) \d(\lambda^{-n} s)$ 
for any polynomial $p(s) \in B = A[s]$. Likewise
$\opn{Tr}_{h_{\lambda} / A} \bigl( t^i \d t \bigr) = 
\lambda^{-(i+1)} t^i \d t$. 
We conclude that
\[ \begin{aligned}
& \lambda^{-(i+1)} p_i(s) \d s = 
(\opn{Tr}_{f' / A} \circ \opn{Tr}_{h_{\lambda} / A})
(t^i \d t) = \\
& \quad 
(\opn{Tr}_{g_{\lambda} / A} \circ \opn{Tr}_{f' / A})(t^i \d t) = 
p_i(\lambda^{-n} s) \d( \lambda^{-n}s) .
\end{aligned} \]
Therefore
$p_i(\lambda^{-n} s) = \lambda^{n -(i+1)} p_i(s)$. 
Since this is true for infinitely many $\lambda$ we must have
$p_i(s) = 0$ for $0 \leq i \leq n-2$, 
and $p_{n-1}(s)$ is a constant. 

In order to compute the value of $p_{n-1}(s) \in A$ we 
note that ${f'}^* \bigl( \d s \bigr) = n t^{n-1} \d t$.
Since we are in characteristic $0$ we can divide by $n$, and 
by Proposition \ref{prop5.1} we get
\[ \begin{aligned}
& \opn{Tr}_{C' / B' / A} \bigl( t^{n-1} \d t \bigr) = 
\opn{Tr}_{C' / B' / A} \bigl({f'}^* \bigl( n^{-1} \d s \bigr) 
\bigr) \\
& \qquad =
\opn{tr}_{C' / B' / A}(1_{C'}) \cdot n^{-1} \d s
= \d s .
\end{aligned} \]
\end{proof}

\begin{rem}
The extra assumptions on the algebra $A$ in the last two 
propositions 
are not really necessary. In Proposition \ref{prop5.2}
we can actually let $A$ be an arbitrary algebra in 
$\cat{EFTAlg} / \K$. For the proof we would then use 
``rigid base change'' (which is developed in \cite{Ye4} to prove 
results on the residue map). Base change allows us to replace both 
$\K$ and $A$ with the ring of integers $\mbb{Z}$. 

Similarly, in Proposition \ref{prop5.1} we can let 
$A$ be an arbitrary algebra in $\cat{EFTAlg} / \K$.
However the proof here requires methods that aren't available 
yet, namely rigid complexes over adically complete rings. 
The idea is 
to use rigid base change, and the setup explained in Remark 
\ref{rem5.1}, to replace $A \to B \to C$ with
$\til{A} \to \til{B} \to \til{C}$. It is possible to choose 
the complete regular local ring 
$\til{A}$ such that it has field of fractions of characteristic 
$0$.
\end{rem}


\end{document}